\documentclass[12pt,letterpaper,english]{article}
\usepackage[latin9]{inputenc}
\usepackage{babel}
\usepackage{amsmath}
\usepackage{amssymb}
\usepackage{graphicx}
\usepackage{setspace}
\usepackage{esint}
\usepackage[authoryear]{natbib}
\usepackage[unicode=true,pdfusetitle,
 bookmarks=true,bookmarksnumbered=false,bookmarksopen=false,
 breaklinks=false,pdfborder={0 0 0},pdfborderstyle={},backref=false,colorlinks=false]
 {hyperref}
\usepackage{breakurl}

\makeatletter

\usepackage{sectsty}
\usepackage{pslatex}
\usepackage{babel}
\usepackage{chngcntr}
\usepackage{apptools}
\usepackage{bm}
\usepackage{amsmath}

\AtAppendix{\counterwithin{lemma}{section}}
\AtAppendix{\counterwithin{proposition}{section}}
\AtAppendix{\counterwithin{theorem}{section}}
\AtAppendix{\counterwithin{corollary}{section}}
\AtAppendix{\counterwithin{assumption}{section}}
\AtAppendix{\counterwithin{definition}{section}}
\AtAppendix{\counterwithin{equation}{section}}

\sectionfont{\sffamily\large}
\subsectionfont{\sffamily\normalsize}
\subsubsectionfont{\sffamily\small}

\newcounter{eqnum}[section] 
\newtheorem{proposition}{Proposition}
\newtheorem{theorem}{Theorem}
\newtheorem{lemma}{Lemma} 
\newtheorem{corollary}{Corollary}
\newtheorem{assumption}{Assumption}

\newcommand{\eproof}{\mbox{}\hfill{\rule{8pt}{8pt}}}

\makeatother

\begin{document}

\title{Game of Variable Contributions to the Common Good under Uncertainty}

\author{H. Dharma Kwon\thanks{Gies College of Business, University of Illinois at Urbana-Champaign,
Champaign, Illinois 61820 \protect \\
Kellogg School of Management, Northwestern University, Evanston, Illinois
60208\protect \\
Email: dhkwon@illinois.edu} }

\date{March 31, 2019}
\maketitle
\begin{abstract}
We consider a stochastic game of contribution to the common good in
which the players have continuous control over the degree of contribution,
and we examine the gradualism arising from the free rider effect.
This game belongs to the class of variable concession games which
generalize wars of attrition. Previously known examples of variable
concession games in the literature yield equilibria characterized
by singular control strategies without any delay of concession. However,
these no-delay equilibria are in contrast to mixed strategy equilibria
of canonical wars of attrition in which each player delays concession
by a randomized time. We find that a variable contribution game with
a single state variable, which extends the Nerlove-Arrow model, possesses
an equilibrium characterized by regular control strategies that result
in a gradual concession. This equilibrium naturally generalizes the
mixed strategy equilibria from the canonical wars of attrition. Stochasticity
of the problem accentuates the qualitative difference between a singular
control solution and a regular control equilibrium solution. We also
find that asymmetry between the players can mitigate the inefficiency
caused by the gradualism. \\
Keywords: Nerlove-Arrow model, war of attrition, stochastic control
game, free rider problem, gradualism
\end{abstract}
\vspace{0.5in}

\newpage{}

\section{Introduction\label{sec:Intro}}

Many business or public policy decisions concern the free rider problem
when contributing to a stock of common good. \citep{Hardin1968}.
It is well-known that a free rider problem induces a wait and see
approach of the individuals who are in a position to contribute to
the common good \citep{Tirole2017}. The wait and see approach in
turn results in underinvestment in the common good. Hence, it is important
for decision makers and social planners to understand the game-theoretic
implications of the free rider problem involving the common good.
Industry examples of a free rider problem with the common good can
be found in the context of generic advertisement for commodities.
For instance, the advertising expenditures by Florida orange juice
advertising programs not only benefit the Florida orange juice industry,
but it also benefits non-Florida orange juice importers \citep{Lee1988}.
In another example, it has been shown that a salmon promotion program
conducted by Norway has benefited its international competitors, too
\citep{Kinnucan2003}. In these examples, the advertising expenditures
of one agent contribute to the stock of the product's overall \emph{goodwill},
``which summarizes the effects of current and past advertising outlays
on demand'' \citep{Nerlove1962}. The stock of goodwill is the common
good in the context of generic advertising because it benefits other
agents, even if they do not contribute to it. In this paper, we examine
the game of variable contribution to the common good where the stock
of common good evolves stochastically. In particular, we obtain the
free rider effect on its Markov perfect equilibrium (MPE) and compare
and contrast it to other games of concession. We address the question
of whether the equilibrium suffers from the gradualism of the players'
contributions to the common good and, if so, whether the inefficiency
arising from the gradualism can be mitigated.

One objective of the paper is to fill the gaps in the equilibrium
characteristics of variable concession games in which the cost is
linear in the contribution. The problem of contribution to common
good belongs to the class of variable concession games in which the
players can control the degree of concession. This class of games
constitutes a significant generalization of the war of attrition.
In the canonical war of attrition, each player can either continue
the game or concede completely, and it typically yields a mixed strategy
equilibrium in which the players delay their concession by a randomized
time. The central question of this paper concerns the characteristics
of the variable concession games. We can imagine three possibilities
of equilibrium strategies: (1) singular control (lump-sum contribution)
strategies without time delay, (2) singular control strategies with
time delay, and (3) regular control strategies that lead to an equilibrium
characterized by gradualism. In the current literature, the variable
concession games thus far have resulted in type (1) equilibrium with
singular control strategies of immediate lump-sum concession. Type
(2), if it exists, is closest to the mixed strategy equilibrium of
the war of attrition, but it has not been found in the literature
or in this paper. The other natural generalization of the mixed strategy
time delay equilibrium to variable concession games is type (3), which
has yet to be found in the current literature on variable concession
games with linear cost. Our paper shows that type (3) is found in
a very simple game-theoretic and stochastic extension of the Nerlove-Arrow
model of goodwill stock \citep{Nerlove1962,Sethi1977,Lon2011}. 

In our model, two players are considering irreversible and costly
contribution to the stock of common good. Each player can contribute
any amount to the common good at any point in time, but the common
good increases the flow profit to both players. The stock of common
good evolves stochastically, and it tends to decline in time on average
unless someone contributes to it, just as the stock of goodwill for
a product depreciates in time without advertisement \citep{Nerlove1962}.
In this game, the strategy of each player is represented by the dynamic
path of its cumulative contribution. We formulate the problem as a
stochastic control game and utilize the well-established stochastic
control theory. In order to find the equilibrium, we need to obtain
the best responses, so we establish the verification theorem for the
best response stochastic control. 

This paper has three main contributions. First, we show that the model
that we consider has a gradualist equilibrium characterized by regular
control. This result is in contrast to the typical control solution:
in a control problem with a linear cost structure, the single decision
maker solution is characterized by singular control rather than regular
control. Second, we find that stochasticity and asymmetry have significant
impact on the equilibrium characteristics. In the deterministic game,
both the singular control solution and the equilibrium solution exhibit
a stable steady state so that an outsider may not be able to tell
the difference between the two. In contrast, in the stochastic case,
the two solutions exhibit markedly different behavior and are easier
to observe in the empirical sense. We also find that asymmetry between
the players destabilizes the gradualist equilibrium, and the outcome
is an asymmetric equilibrium with singular control strategy adopted
by at least one of the players. Hence, asymmetry can mitigate the
inefficiency of the gradualist equilibrium. Third, the paper provides
a mathematical framework to obtain an MPE of a stochastic game of
variable concession involving both singular and regular control. 

Although there are many equilibrium solution concepts, we limit our
attention to MPE \citep{Maskin2001}. MPE is a subgame perfect equilibrium
in which the players' actions are determined by the current value,
but not by the past history, of the economically relevant state variable,
and hence it is a key notion for analyzing a game. 

Cooperative equilibrium concepts are beyond the scope of this paper.
Coordinated plans of action do produce an efficient outcome, which
will change the form of the solution; for instance, the singular control
boundary will change. Although cooperation does happen between contributors
to common good, it often requires prior coordination or bargaining,
and we can still consider the non-cooperative MPE a baseline solution
prior to coordination. For instance, in a Nash bargaining solution
\citep{Nash1950}, the non-cooperative MPE outcome can serve as the
disagreement point, and therefore, it is still a meaningful reference
point.

The paper contributes to the literature on variable concession games,
which an extension of a war of attrition \citep{Smith1974}. Typical
attrition games under complete information possess mixed strategy
equilibria with random time delays, both in the deterministic case
\citep{Hendricks1988} and in the stochastic case \citep{Steg2015,Georgiadis2019}.
In contrast, the known examples of the game of variable concession
exhibit singular control equilibria with no time delay. One example
is Cournot competition under declining demand when the firms can reduce
the production capacity at a variable cost \citep{Ghemawat1990}.
The equilibrium strategy is to immediately reduce the capacity to
the myopic Cournot equilibrium level through singular control. The
stochastic generalization of the Cournot model also exhibits similar
characteristics \citep{Steg2012}.

The paper also contributes to the literature on games of contribution
to public goods. \citet{Fershtman1991} examine a dynamic game of
voluntary contribution to public goods. In their model, players continuously
contribute to the stock of public goods over time. They obtain an
equilibrium using a differential game approach and demonstrate that
the free riding problem is acute without commitment. Their model is
similar to ours, but they model a situation with costs that grow quadratically
with the rate of contribution, so it is prohibitively costly to make
a lump sum contribution. Since our model allows for a lump sum contribution
due to the linear cost structure, the characteristics of the equilibria
are very different, and it is difficult to compare their results to
ours. \citet{Battaglini2014} examine a problem of dynamic free riding
in which each individual allocates its endowment between private consumption
and irreversible contribution to the public good. They study the implications
of the irreversibility of their model and conclude that irreversibility
can alleviate inefficiency of the equilibria. It is noteworthy that
the equilibrium of their model involves lump sum contribution (singular
control) strategies. \citet{Ferrari2017} examine a significantly
generalized model with stochasticity to obtain its equilibria and
study the effect of uncertainty and irreversibility of contribution
to the public good. They also obtain equilibria characterized by lump
sum contribution strategies. In contrast to the literature on private
consumption and contribution to the public good, our model does not
incorporate consumption of the players.

Because our model assumes that the cost of contribution is linear
in the magnitude of the improvement in the common good, we formulate
it as a game-theoretic extension of monotone follower singular control
problems with a single dimensional state variable. In a similar vein,
\citet{Lon2011} apply singular control framework to the Nerlove-Arrow
model of expenditure in the stock of goodwill. Recently, some work
on game-theoretic study of singular control problems has emerged.
\citet{Steg2012} examines Cournot competition that leads to a singular
control equilibrium. \citet{Kwon2015} examines a singular control
game in the context of a market share competition in which a player's
control is to negate his opponent's payoff. \citet{Ferrari2017} also
analyze a model that incorporates game-theoretic singular control,
but the model has a more complex structure as the players make consumption
and contribution decisions at the same time. 

The paper is organized as follows. In Section \ref{sec:VariableContribution},
we examine a game of variable contribution to the common good and
show that it yields a regular control strategy equilibrium. In Section
\ref{sec:Impact-of-Stoch}, we examine the impact of stochasticity
and asymmetry between the players. In particular, we show that asymmetry
eliminates the regular control equilibrium thereby improving the efficiency.
In Section \ref{sec:Discussions}, we discuss several aspects of the
results that are worthy of note. In Section \ref{sec:Conclusions},
we summarize the main results and implications of the paper and provide
concluding remarks. 

\section{Variable Contribution Game\label{sec:VariableContribution}}

In this section, we present a game of variable contribution to the
common good that results in a regular control equilibrium. We first
present the model in Section \ref{subsec:The-Model}, and then we
examine the single decision maker case as a benchmark in Section \ref{subsec:Single}.
We construct the verification theorem for best responses in Section
\ref{subsec:Game} and obtain the regular control MPE in Section \ref{subsec:Symmetric-Eq}. 

\subsection{The Model \label{subsec:The-Model}}

We consider a game between two players, each of whom receives a flow
profit that depends on a \emph{common} state variable. Either player
can boost the common state variable at a cost by any amount at any
point in time. The model is applicable to a number of industry examples.
One example is a game between two manufacturers who share a common
supplier. Each manufacturer can make a variable investment to boost
the quality of the shared supplier, which in turn benefits the other
manufacturer through spillover \citep{Muthulingam2016,Kim2017a}.
Another example is the game of irreversible and variable investment
in the stock of goodwill \citep{Nerlove1962} through advertisement
such as in generic advertising on commodities \citep{Lee1988,Kinnucan2003}. 

We let the process $Z=\{Z_{t}:t\ge0\}$ denote the stock of common
good defined in the interval $I=(a,b)\subseteq\mathbb{R}$ on a filtered
probability space $(\Omega,\mathcal{F},\mathcal{F}_{t},\mathbb{P})$
that satisfies the \emph{usual condition}. If $I=\mathbb{R}$, for
example, it is understood that $a=-\infty$ and $b=\infty$. We assume
that $Z$ satisfies the following stochastic differential equation
(SDE): 
\[
dZ_{t}=\mu(Z_{t})dt+\sigma(Z_{t})dW_{t}+d\xi_{it}+d\xi_{jt}\;,
\]
where $W=\{W_{t}:t\ge0\}$ is a Wiener process progressively measurable
with respect to $\{\mathcal{F}_{t}:t\ge0\}$. Here $\mu(\cdot)$ is
the drift term which we interpret as the time-averaged rate of change
of $Z$ in the absence of control. In this paper, we assume $\mu(\cdot)<0$
to model the deterioration of the common good. The volatility $\sigma(\cdot)>0$
represents the magnitude of the white noise. The process $\xi_{i}=\{\xi_{it}:t\ge0\}$
is a non-decreasing càdlàg (right continuous with left limits) process
controlled by player $i$ adapted to $\{\mathcal{F}_{t}:t\ge0\}$.
We interpret $\xi_{it}$ as the cumulative contribution of player
$i$ to $Z$ up to time $t$. Since each player $i$ controls the
process $\xi_{i}$, we say that $\xi_{i}$ is player $i$'s strategy,
and $\xi=(\xi_{i},\xi_{j})$ is the strategy profile. Throughout the
paper, we let $\Xi_{i}$ denote the set of all possible $\mathcal{F}_{t}$-adapted
control processes $\xi_{i}$.

We remark that $\xi_{i}$ is composed of a continuous process and
a discontinuous process as follows:
\[
\xi_{it_{2}}-\xi_{it_{1}^{-}}=\int_{t_{1}}^{t_{2}}d\xi_{it}^{c}+\sum_{t\in[t_{1},t_{2}]}\Delta\xi_{it}\:,
\]
where $\xi_{i}^{c}$ is the continuous part of $\xi_{i}$, and $\Delta\xi_{it}=\xi_{it}-\xi_{it^{-}}$
is the instantaneous jump in $\xi_{i}$ at time $t$. Similarly, we
can decompose the process $Z_{t}$ into a continuous part $Z_{t}^{c}$
and a discontinuous part $\Delta Z_{t}=Z_{t}-Z_{t^{-}}=\Delta\xi_{1t}+\Delta\xi_{2t}$. 

Given a strategy profile $\xi$, player $i$'s payoff is given by
the following function:
\[
V_{i}(x;\xi)=\mathbb{E}^{x}\left[\int_{0}^{\infty}e^{-rt}\pi_{i}(Z_{t})dt-\int_{0}^{\infty}e^{-rt}k_{i}d\xi_{it}\right]\:.
\]
 Here $\mathbb{E}^{x}[\cdot]=\mathbb{E}[\cdot\vert Z_{0}=x]$ is the
conditional expectation operator given the initial condition $Z_{0}=x$.
The integrand $\pi_{i}(\cdot)$ is a non-decreasing function that
represents the profit flow for player $i$, and $k_{i}>0$ is the
cost of increasing a unit of $\xi_{i}$. Lastly, $r>0$ is the discount
rate common to both players. 

For the sake of analytical tractability, we make a number of assumptions
below that are standard in the stochastic control literature. We first
make some assumptions regarding $\mu(\cdot)$ and $\sigma(\cdot)$.
Let $X=\{X_{t}\,:\,t\ge0\}$ denote the uncontrolled process which
satisfies $dX_{t}=\mu(X_{t})dt+\sigma(X_{t})dW_{t}$.

\begin{assumption} \label{assump:mu-sigma} (i) $\mu(\cdot)$ and
$\sigma(\cdot)$ are Lipschitz continuous functions satisfying $\vert\mu(x)\vert+\vert\sigma(x)\vert\le\delta(1+\vert x\vert)$
for some constant $\delta>0$. 

(ii) $\{e^{-r\tau}(X_{t})^{-}:\text{\ensuremath{\tau\ }}\text{is a stopping time},\;\tau<\infty\}$
is uniformly integrable for any initial value $X_{0}=x$. Furthermore,
$\lim_{t\rightarrow\infty}\mathbb{E}^{x}[e^{-rt}(X_{t})^{-}]=0$.
\end{assumption}

Assumption \ref{assump:mu-sigma} (i) implies that the uncontrolled
process $X$ has a unique strong solution to the SDE. It also implies
that $\sigma(\cdot)$ is locally bounded; this will be useful when
we apply Dynkin's formula to the payoff function because the stochastic
integral involving $\sigma(Z_{t})dW_{t}$ is a \emph{local martingale}
which possesses convenient properties (Chapter IV, \citealp{Revuz1999}).
Assumption \ref{assump:mu-sigma} (ii) ensures that the limiting behaviors
of the process $X$ are well-defined so that we can construct a verification
theorem in Section \ref{subsec:Game}. 

Below we let $C(I)$ denote the set of continuous functions defined
on $I$.

\begin{assumption} \label{assump:abs-int} $\pi_{i}(\cdot)\in C(I)$
is strictly increasing and bounded from above, i.e., $\lim_{x\rightarrow b}\pi_{i}(x)<\pi_{M}$
for some positive constant $\pi_{M}$. Furthermore, it satisfies the
absolute integrability condition $\mathbb{E}^{x}\left[\int_{0}^{\infty}\left|e^{-rt}\pi_{i}(X_{t})\right|dt\right]<\infty$
for the uncontrolled process $X$. \end{assumption} Assumption \ref{assump:abs-int}
ensures that the payoff $V_{i}(x;\xi)$ is well-defined and that the
function
\begin{equation}
(R\pi_{i})(x):=\mathbb{E}^{x}\left[\int_{0}^{\infty}e^{-rt}\pi_{i}(X_{t}^{x})dt\right]\label{eq:Rpi}
\end{equation}
exists. The function $(R\pi_{i})(\cdot)$ has the meaning of the payoff
from perpetually keeping an uncontrolled process $X$. Later we establish
that $(R\pi_{i})(\cdot)$ is an element of the payoff function.

Next, we define the $r$-excessive characteristic operator \citep{Alvarez2003}:
\begin{equation}
\mathcal{A}:=\frac{1}{2}\sigma(x)^{2}\partial_{x}^{2}+\mu(x)\partial_{x}-r\:.\label{eq:char-op-A}
\end{equation}
 We let $\psi(\cdot)$ and $\phi(\cdot)$ respectively denote two
linearly independent increasing and decreasing fundamental solutions
to the differential equation $\mathcal{A}\psi(x)=\mathcal{A}\phi(x)=0$
\citep{Borodin1996,Alvarez2003}. We remark that $(R\pi_{i})(\cdot)$
satisfies the differential equation $\mathcal{A}(R\pi_{i})(x)+\pi_{i}(x)=0$
according to the optimal stopping theory. 

We also define the following functions
\begin{align}
q_{i}(x) & =\pi_{i}(x)+\mathcal{A}k_{i}x=\pi_{i}(x)+[\mu(x)-rx]k_{i}\;,\label{eq:q1-x-1}\\
(Rq_{i})(x) & =\mathbb{E}^{x}\left[\int_{0}^{\infty}e^{-rt}q_{i}(X_{t})dt\right]\:,\nonumber 
\end{align}
 and assume the following properties of $(R\pi_{i})(\cdot)$ and $q_{i}(\cdot)$:

\begin{assumption} \label{assump:qi} (i) $(R\pi_{i})(\cdot)$ satisfies
$\lim_{x\rightarrow a}(R\pi_{i})(x)<0$ and $\lim_{x\rightarrow a}(R\pi_{i})(x)/\phi(x)=0$.
(ii) There exists some $x_{i}^{*}\in I$ such that $q_{i}(x)$ is
strictly increasing for $x<x^{*}$ and strictly decreasing for $x>x^{*}$.
(iii) $\lim_{x\rightarrow a}q(x)=-\infty$, $\lim_{x\rightarrow a}q(x)\psi'(x)\exp\left(\int_{0}^{x}\frac{2\mu(y)}{\sigma^{2}(y)}dy\right)=0$,
$\lim_{x\rightarrow b}q(x)\phi'(x)\exp\left(\int_{0}^{x}\frac{2\mu(y)}{\sigma^{2}(y)}dy\right)=0$,
and $\lim_{x\rightarrow b}(Rq)'(x)/\phi'(x)>0$. \end{assumption}

\noindent Assumption \ref{assump:qi} serves as the sufficient condition
for the unique optimal control solution to exist for the model examined
in Section \ref{subsec:Single}. Specifically, in the single decision
maker's problem, the assumptions drive a solution with a singular
control region of the form $(a,\theta)$ for some threshold $\theta$.
Assumption \ref{assump:qi} (i) ensures that the flow profit function
$\pi_{i}(\cdot)$ is negative and well-behaved near $a$. To gain
an intuitive understanding of the assumptions regarding $q_{i}(\cdot)$,
we consider a special case of $\mu(x)=\mu$, in which case $q_{i}'(x)=\pi_{i}'(x)-rk_{i}$.
Then Assumption \ref{assump:qi} (ii) implies that $\pi_{i}'(x)>rk_{i}$
for $x<x_{i}^{*}$ and $\pi_{i}'(x)<rk_{i}$ for $x>x_{i}^{*}$. This
implies that it is optimal to boost $\xi_{i}$ if and only if $x<x_{i}^{*}$.
Thus, the players have incentive to boost $\xi_{i}$ only if $X_{t}$
falls below a threshold. Assumption \ref{assump:qi} (iii) ensures
that a unique threshold for boosting $\xi_{i}$ exists for the model
examined in Section \ref{subsec:Single}. 

\subsection{Benchmark: Single Decision Maker Problem \label{subsec:Single}}

In this subsection, we review the single decision maker problem as
a benchmark and provide the optimal solution. This class of problems
is extensively examined in the literature \citep{Alvarez2001,Oksendal2005,Lon2011},
but we reproduce it here because its solution will be utilized in
the equilibrium solution of the game in the remainder of the paper. 

Since there is only one decision maker, we drop the player index $i$
for convenience throughout Section \ref{subsec:Single}. The value
function associated with a control policy $\xi=\{\xi_{t}:t\ge0\}$
is given by 
\[
V(x;\xi)=\mathbb{E}^{x}\left[\int_{0}^{\infty}e^{-rt}\pi(Z_{t})dt-\int_{0}^{\infty}e^{-rt}kd\xi_{t}\right]\:.
\]
 The objective of the decision maker is to maximize $V(\cdot;\xi)$
with respect to $\xi$. This class of problems is known as the singular
stochastic control monotone follower problems. 

We first provide the sufficient condition (the \emph{optimality condition})
for the optimal control solution. We let $C^{n}(I)$ denote the set
of functions defined on $I$ that are $n$ times continuously differentiable.
Suppose that $V(\cdot)\in C^{2}(I)$ satisfies the following conditions:
(i) $\mathcal{A}V(x)+\pi(x)\le0$ and $V'(x)-k\le0$ for all $x\in I$,
and (ii) $[\mathcal{A}V(x)+\pi(x)][V'(x)-k]=0$ for all $x\in I$.
Then $V(x)$ coincides with the optimal solution $\sup_{\xi\in\Xi}V(x;\xi)$.
The proof of this sufficient condition is provided, for example, by
\citet{Oksendal2005} and \citet{Lon2011}, and so we will not reproduce
it here. 

\begin{lemma} \label{lemma:OptSol} Under Assumptions \ref{assump:mu-sigma}\textendash \ref{assump:qi},
there exist a threshold $\theta\in I$ and a coefficient $A$ such
that the optimal solution $V(\cdot)\in C^{2}(I)$ is given as follows: 

\begin{equation}
V(x)=\begin{cases}
k(x-\theta)+V(\theta) & \text{for}\;x<\theta\\
(R\pi)(x)+A\phi(x) & \text{for}\;x\ge\theta
\end{cases}\:.\label{eq:Vx}
\end{equation}
\end{lemma}

Next, we provide the intuition behind the solution through a numerical
example. 

\emph{Example 1}: We consider a problem of constant $\mu(x)=\mu<0$
and $\sigma(x)=\sigma>0$ and the flow profit $\pi(x)=1-\exp(\nu x)$
with $\nu<0$ that satisfies $\frac{1}{2}\sigma^{2}\nu^{2}+\mu\nu-r<0$.
Then it is straightforward to verify that $\phi(x)=\exp(\gamma_{-}x)$
and $\psi(x)=\exp(\gamma_{+}x)$ where
\[
\gamma_{\pm}=\frac{1}{\sigma^{2}}(-\mu\pm\sqrt{\mu^{2}+2r\sigma^{2}})\:,
\]
and 
\begin{align*}
(R\pi)(x) & =\frac{1}{r}-\frac{\exp(\nu x)}{\beta}\:,\\
\text{where }\quad\beta & =r-\mu\nu-\frac{1}{2}\sigma^{2}\nu^{2}\:.
\end{align*}
 Assumptions \ref{assump:mu-sigma}\textendash \ref{assump:qi} are
satisfied in this example so that Lemma \ref{lemma:OptSol} applies.
Furthermore, the optimal threshold $\theta$ and the coefficient are
given by 
\begin{align*}
\theta & =\frac{1}{\nu}\ln\frac{k\beta\gamma_{-}}{\nu(\nu-\gamma_{-})}\;,\\
A & =\frac{k-(R\pi)'(\theta)}{\phi'(\theta)}\:.
\end{align*}

\begin{figure}
\begin{centering}
\includegraphics[scale=0.3]{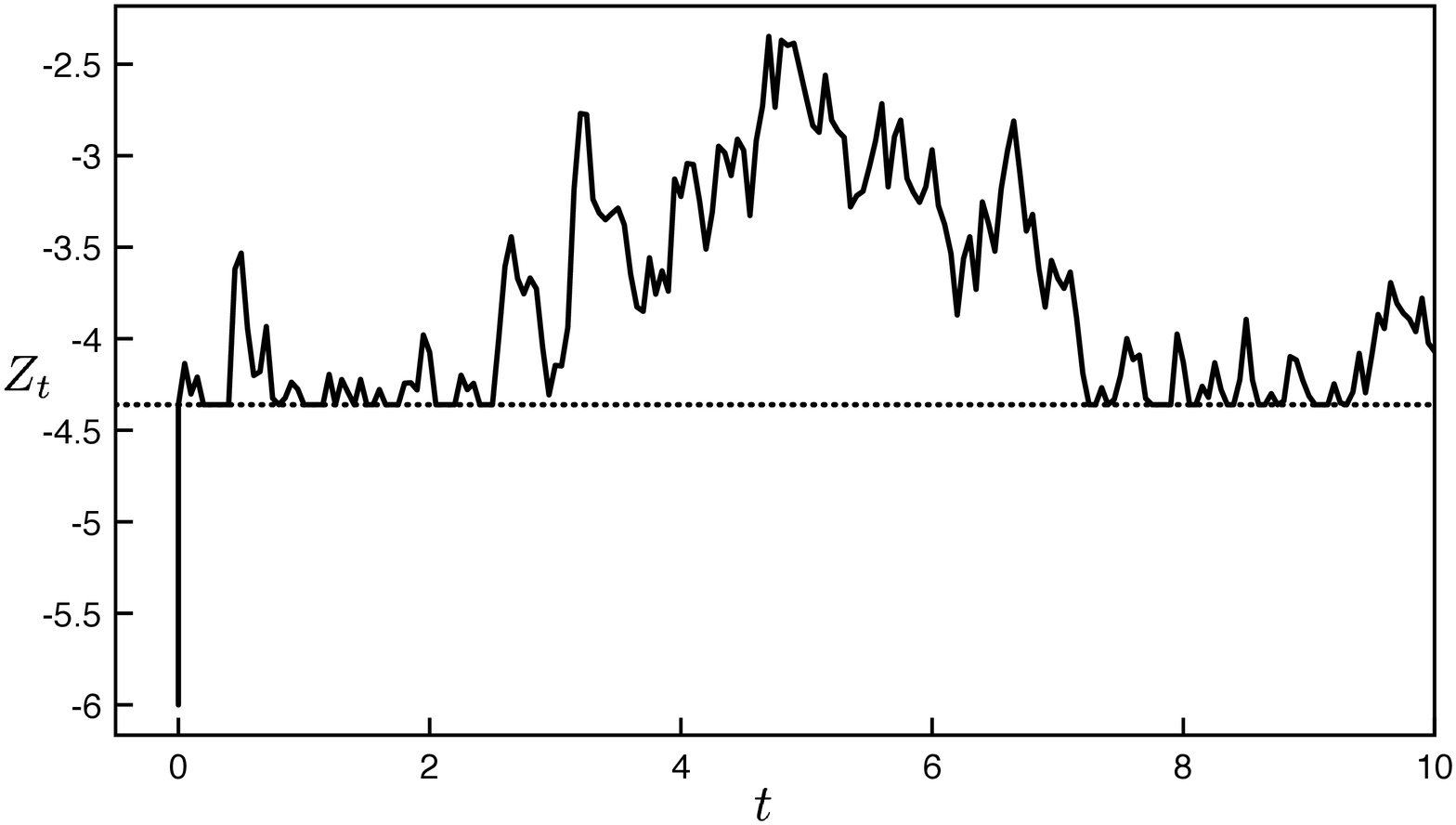}
\par\end{centering}
\begin{centering}
\includegraphics[scale=0.3]{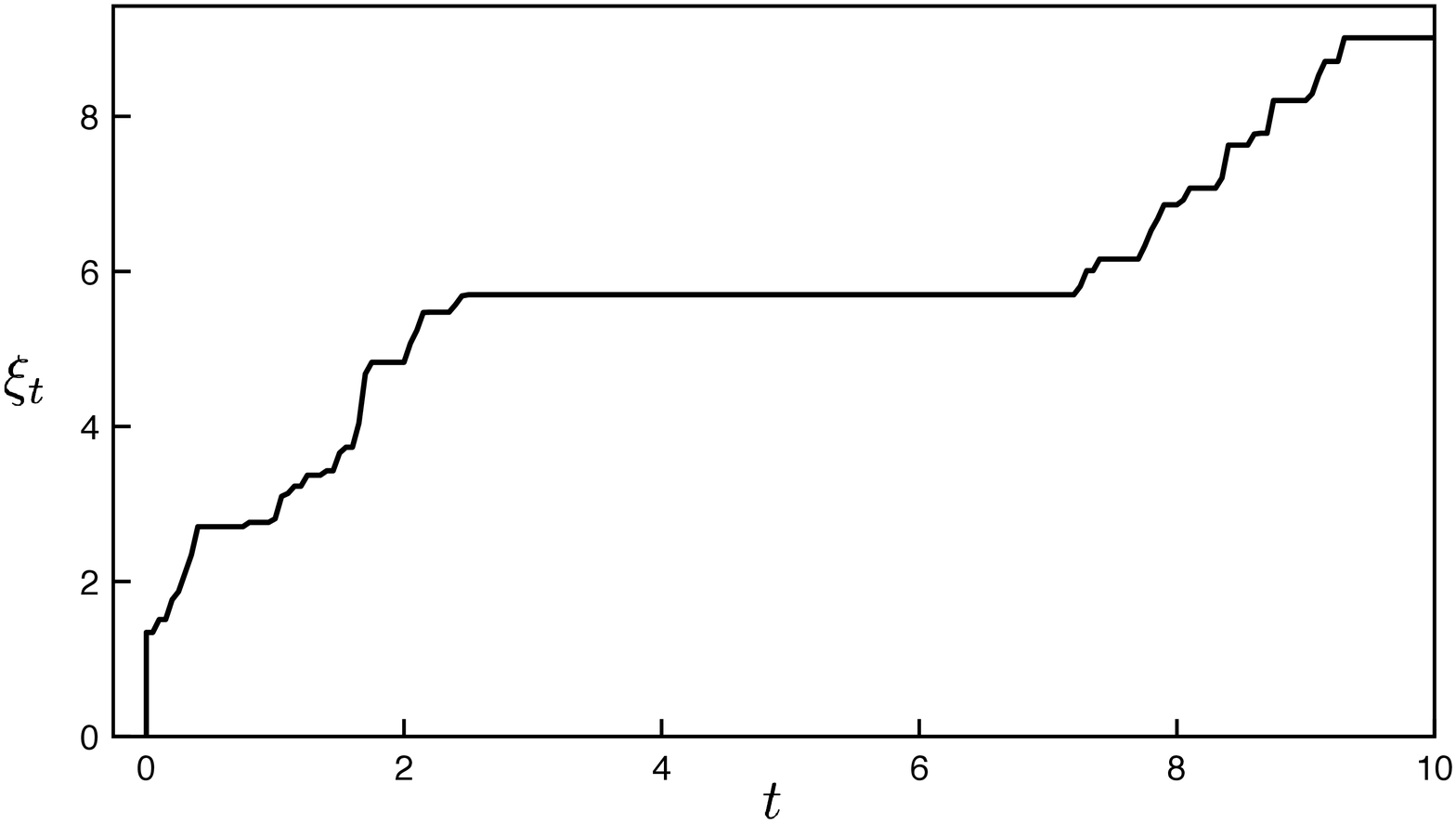}
\par\end{centering}
\begin{centering}
\includegraphics[scale=0.3]{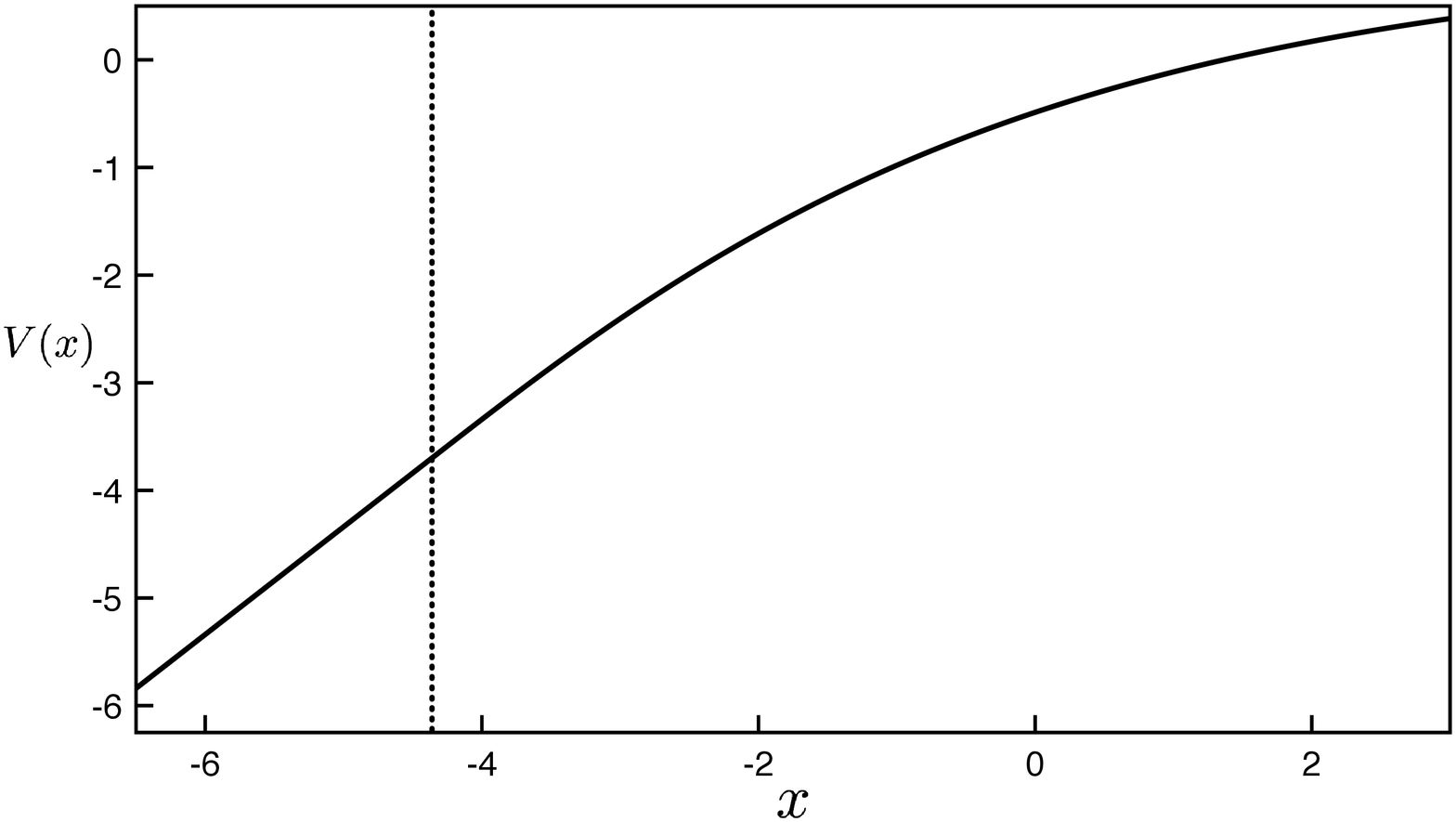}
\par\end{centering}
\caption{Single decision maker solution. The dotted lines indicate the threshold
$\theta$. \label{fig:Single-decision-maker}}
\end{figure}

For the numerical illustration in Fig. \ref{fig:Single-decision-maker},
we set $\mu=-1$, $\sigma=r=k=1$. In this case, $\theta=-4.3605$
is the optimal threshold. See Fig. \ref{fig:Single-decision-maker}
for a simulated sample path of $Z$ and $\xi$ and the optimal value
function $V(x)$. 

The optimal policy associated with the value function (\ref{eq:Vx})
is a singular control policy: to boost $Z_{t}$ instantaneously up
to $\theta$ whenever $Z_{t}$ falls below $\theta$. If the initial
value $Z_{0}$ is below $\theta$, then $Z_{t}$ discontinuously jumps
to $\theta$ at time $0$, after which $Z_{t}$ is continuous. The
threshold $\theta$ functions as a reflecting boundary as shown in
Fig. \ref{fig:Single-decision-maker}. In the region below $\theta$,
the value function is linear in $x$ with the slope $k$, as illustrated
in the figure as well as expressed in (\ref{eq:Vx}). This is because
$(a,\theta)$ is the singular control region in which it costs exactly
$k$ to boost $Z$ by a unit. 

\subsection{Verification Theorem for Best Responses \label{subsec:Game}}

Next, we return to the game-theoretic model introduced in the beginning
of the section. Our goal is to construct the verification theorem
for best responses, which will then be used to construct MPEs in the
remainder of the paper. 

In a conventional solution to the singular control problem as in Section
\ref{subsec:Single}, the optimal control process is decomposed as
$d\xi_{it}=d\xi_{it}^{l}+\Delta\xi_{it}$ where $\Delta\xi_{it}=\xi_{it}-\xi_{it^{-}}$
represents the discontinuous evolution of the process $\xi_{i}$ while
$\xi_{it}^{l}$ is a continuous process like a \emph{local time} \citep{Protter2003}.
The process $\xi_{it}^{l}$ is not absolutely continuous with respect
to Lebesgue measure and cannot be represented as an integral $\int_{0}^{t}u_{is}ds$
for any process $u_{i}=\{u_{it}:t\ge0\}$ \citep{Karatzas1983}; for
example, the sample path of $\xi$ in Figure \ref{fig:Single-decision-maker}
does not possess well-defined time derivatives when $\xi_{t}$ increases
in time. Both components $\xi_{it}^{l}$ and $\Delta\xi_{it}$ constitute
the singular part of $\xi_{it}$. In general, however, a control process
must also encompass a regular control process as follows:
\[
d\xi_{it}=u_{it}dt+d\xi_{it}^{l}+\Delta\xi_{it}\;,
\]
where $u_{it}\ge0$ is a process adapted to $(\mathcal{F}_{t})$. 

To apply the conventional stochastic control theory, we need to define
the feasible space of $u_{i}$. More specifically, in order for the
SDE $dZ_{t}=\mu(Z_{t})dt+\sigma(Z_{t})dW_{t}+\sum_{i=1}^{2}d\xi_{it}$
to have a unique strong solution, we need to limit $u_{i}$ within
the class of $u_{i}(t,x)$ that satisfies the two following conditions:
(1) $u_{i}(t,x)$ is Lipschitz continuous in $x$, and (2) $\vert u_{i}(t,x)\vert\le\delta(1+\vert x\vert)$
for some constant $\delta$. Let $\mathcal{U}$ be the set of functions
$u:(0,\infty)\times I\rightarrow\mathbb{R}$ that satisfy these two
conditions. Then we let $\Sigma_{i}$ denote the set of $\mathcal{F}_{t}$-adapted
processes $\xi_{it}$ that satisfy $d\xi_{it}=u_{i}(t,Z_{t})dt+d\xi_{it}^{l}+\Delta\xi_{it}$
for some $u_{i}\in\mathcal{U}$. We remark that if $\text{\ensuremath{\xi}}_{i}\in\Sigma_{i}$
then the SDE of $Z$
\[
dZ_{t}=[\mu(Z_{t})+\sum_{i=1}^{2}u_{i}(t,Z_{t})]dt+\sigma(Z_{t})dW_{t}+\sum_{i=1}^{2}(d\xi_{it}^{l}+\Delta\xi_{it})
\]
satisfies the sufficient condition for possessing a unique strong
solution because $u_{i}\in\mathcal{U}$. Note also that $\Sigma_{i}$
is a proper subset of the feasible strategy space $\Xi_{i}$, which
is the set of all possible $\mathcal{F}_{t}$-adapted control processes
$\xi_{it}$. 

Let $\mathcal{M}_{i}$ be the set of player $i$'s Markov control
strategies $\xi_{i}$ which depend only on the current value of $Z_{t}$.
It means that $\xi_{it}$ satisfies $d\xi_{it}=u_{i}(Z_{t})dt+d\xi_{it}^{l}+\Delta\xi_{it}$
where the singular control region is given as a subset of $I$. For
instance, if the singular control region of player $i$ is $[\alpha,\beta]$,
then whenever $Z_{t^{-}}\in[\alpha,\beta]$, $\xi$ undergoes a jump
$\Delta\xi_{it}=\beta-Z_{t^{-}}$, i.e., player $i$ boosts $Z$ up
to $\beta$. Furthermore, $d\xi_{it}^{l}>0$ only when $Z_{t}$ hits
$\beta$. By definition, an MPE is a subgame perfect equilibrium $\xi=(\xi_{i},\xi_{j})$
that belongs to $\mathcal{M}_{i}\times\mathcal{M}_{j}$. 

For the purpose of obtaining MPE $\xi\in\mathcal{M}_{i}\times\mathcal{M}_{j}$
we can safely focus on obtaining the best response strategies $\xi_{i}'\in\Xi_{i}$
in response to the opponent's Markov strategy $\xi_{j}\in\mathcal{M}_{j}$,
i.e., to obtain $\xi_{i}\in\Xi_{i}$ such that $V_{i}(x;\xi_{i},\xi_{j})=\sup_{\xi_{i}'\in\Xi_{i}}V_{i}(x;\xi_{i}',\xi_{j})$.
If the best response $\xi_{i}$ happens to belong to $\mathcal{M}_{i}$,
then $\xi=(\xi_{i},\xi_{j})$ is an MPE. Note that we do not, however,
look for $\xi_{i}\in\mathcal{M}_{i}$ such that $V_{i}(x;\xi_{i},\xi_{j})=\sup_{\xi_{i}'\in\mathcal{M}_{i}}V_{i}(x;\xi_{i}',\xi_{j})$;
because $\mathcal{M}_{i}\subset\Xi_{i}$, there may be another $\xi_{i}'\in\Xi_{i}$
such that $V_{i}(x;\xi_{i}',\xi_{j})>V_{i}(x;\xi_{i},\xi_{j})$, in
which case $(\xi_{i},\xi_{j})$ is not a Nash equilibrium. 

Since we assume a Markov control process $\xi_{j}\in\mathcal{M}_{j}$,
we can partition the interval $I$ into regions of discontinuous $\xi_{jt}$
and continuous $\xi_{jt}$. Let $C_{j}$ denote the open subset of
$I$ in which $\xi_{jt}$ evolves continuously and non-singularly
in time, and let $D_{j}=I\backslash C_{j}$ denote the \emph{singular
control region} where $d\xi_{jt}^{c}>0$ or $\Delta\xi_{jt}>0$. 

Now we provide sufficient conditions for the best response to $\xi_{j}$. 

\begin{theorem}\label{thm:verif} Suppose Assumptions \ref{assump:mu-sigma}\textendash \ref{assump:qi}
hold. Assume player $j$'s Markov strategy $\xi_{j}\in\mathcal{M}_{j}\cap\Sigma_{j}$
that satisfies $d\xi_{jt}=u_{j}(Z_{t})dt+d\xi_{jt}^{l}+\Delta\xi_{jt}$
for some function $u_{j}(\cdot)\ge0$. Suppose that there exist a
function $U(\cdot)$ on $I$ and some Markov strategy $\xi_{i}^{*}\in\mathcal{M}_{i}\cap\Sigma_{i}$
that satisfy the SDE 
\begin{equation}
d\xi_{it}^{*}=u_{i}^{*}(Z_{t})dt+d\xi_{it}^{*l}+\Delta\xi_{it}^{*}\:\label{eq:xi*}
\end{equation}
for some $u_{i}^{*}(\cdot)\ge0$ and the following conditions:

(i) $U(\cdot)\in C^{2}(C_{j})\cap C(I)$, and $U(\cdot)$ is non-decreasing
and bounded from above. Furthermore, $U'(x)=0$ in the interior of
$D_{j}$, and $U'(Z_{t})d\xi_{jt}^{l}=0$ for all $t$, where $Z$
is the state process that evolves under the strategy profile $(\xi_{i}^{*},\xi_{j})$.

(ii) There is a function $\tilde{U}(x)\in C^{2}(I)$ such that $\tilde{U}(x)=U(x)$
for all $x\in C_{j}$ and $\tilde{U}'(x)$ is bounded for $x\in D_{j}$. 

(iii) $\max\{\mathcal{A}U(x)+\pi_{i}(x)+u_{j}(x)U'(x)+v_{i}U'(x)-v_{i}k_{i},U'(x)-k_{i}\}\le0$
for all $x\in C_{j}$ and any arbitrary $v_{i}\ge0$. 

(iv) Let $D_{i}=\{x\in I:\Delta\xi_{it}^{*}>0\;\text{whenever}\:x=Z_{t^{-}}\}\cap C_{j}$
be player $i$'s singular control region within $C_{j}$ and $C_{i}=\{x\in I:\Delta\xi_{it}^{*}=0\;\text{whenever}\:x=Z_{t^{-}}\}\cap C_{j}$.
Then $\mathcal{A}U(x)+\pi_{i}(x)+u_{j}(x)U'(x)+u_{i}^{*}(x)U'(x)-u_{i}^{*}(x)k_{i}=0$
for all $x\in C_{i}$ and $U'(x)=k_{i}$ for all $x\in D_{i}$. 

Then $\xi_{i}^{*}$ is the best response to $\xi_{j}$ amongst all
control processes that belong to $\Sigma_{i}$, i.e., $V_{i}(x;\xi_{i}^{*},\xi_{j})=\sup_{\xi_{i}\in\Sigma_{i}}V_{i}(x;\xi_{i},\xi_{j})$.
\end{theorem}

\emph{Remark}: Strictly speaking, Theorem \ref{thm:verif} does not
give sufficient conditions for the best response among the whole strategy
space $\Xi_{i}$, but it gives sufficient conditions for the best
response among the limited space $\Sigma_{i}$. However, the strategy
profiles that we obtain in Theorem \ref{thm:symm-eq} and Proposition
\ref{prop:asymm-eq} are proper MPE. This result is obtained because
even though Theorem \ref{thm:verif} is used to obtain $(\xi_{1},\xi_{2})\in\Sigma_{1}\times\Sigma_{2}$
that satisfies $V_{i}(x;\xi_{i},\xi_{j})=\sup_{\xi_{i}'\in\Sigma_{i}}V_{i}(x;\xi_{i}',\xi_{j})$
for both $i=1,2$, we can show that 
\[
V_{i}(x;\xi_{i},\xi_{j})=\sup_{\xi_{i}'\in\Xi_{i}}V_{i}(x;\xi_{i}',\xi_{j})\;,
\]
by the optimality condition for singular control given in Section
\ref{subsec:Single}. The same is true if $i$ and $j$ are interchanged.
Therefore, $(\xi_{1},\xi_{2})$ is a proper MPE. The detail is provided
in the proof of Theorem \ref{thm:symm-eq}. 

\subsection{Regular Control Strategy Equilibrium \label{subsec:Symmetric-Eq}}

Next, we construct a regular control MPE. We define a \emph{regular
control MPE} as one with control strategies of the following form
for both players $i=1,2$:
\[
d\xi_{it}=u_{i}(Z_{t})dt\;,
\]
where $\xi_{i}\in\Sigma_{i}$. In this subsection, we assume that
the two players are symmetric, i.e., $\pi_{i}=\pi_{j}=\pi$ and $k_{i}=k_{j}=k$. 

\begin{theorem} \label{thm:symm-eq} Suppose Assumptions \ref{assump:mu-sigma}\textendash \ref{assump:qi}
hold. Let $V(\cdot)$ denote the solution (\ref{eq:Vx}) to the single
decision maker problem. Also define a symmetric strategy profile $\xi=(\xi_{1},\xi_{2})$
with the regular control process $d\xi_{it}=u(Z_{t})dt$, where 
\begin{equation}
u(x)=\begin{cases}
-\frac{1}{k}\left(\mathcal{A}V(x)+\pi(x)\right) & \text{for}\;x<\theta\\
0 & \text{for}\;x\ge\theta
\end{cases}\:.\label{eq:u-x-reg}
\end{equation}
Lastly, suppose $u(\cdot)\in\mathcal{U}$. Then $\xi$ is an MPE with
a symmetric payoff $V_{i}(x;\xi)=V_{j}(x;\xi)=V(x)$. \end{theorem}

Theorem \ref{thm:symm-eq} obtains an MPE completely characterized
by regular control. Intuitively, both players exert gradual control
$u(\cdot)$ if and only if $Z$ is sufficiently low (less than $\theta$).
The regular control MPE is reminiscent of the mixed strategy delay
equilibrium of the canonical war of attrition in which both players
gradually concede in the probabilistic sense through a Poisson process.
Thus, we may consider the regular control strategy equilibrium a generalization
of the mixed strategy delay equilibrium.

One notable characteristic of the regular control MPE is that the
threshold of the control region $\theta$ is identical to the threshold
of singular control region in the single decision maker solution.
It implies that the free rider effect does not shift the control threshold;
instead, it drives gradualism. 

\emph{Example }2: We now consider the game-theoretic extension of
Example 1. For analytical tractability, the profit flow $\pi(\cdot)$
is modified as follows:
\[
\pi(x)=\begin{cases}
1-\exp(\nu x) & \text{for}\;x\ge x_{c}\\
\pi(x_{c})+(x-x_{c})\rho & \text{for}\;x<x_{c}
\end{cases}\:,
\]
where $\rho>rk$ and $x_{c}<\theta$ are parameters that we specify
below. The form of $\pi(\cdot)$ is modified so that $\vert\pi(x)\vert$
does not grow faster than $\vert x\vert$ for sufficiently large $\vert x\vert$.
The modification is necessary to ensure that $u_{i}(x)$ in equilibrium
does not grow faster than $\vert x\vert$ for large $\vert x\vert$,
preserving the well-known sufficient conditions for existence and
uniqueness of the strong solution to SDE of $Z$. 

In this case, the function $(R\pi)(x)$ has a more complicated form.
Define 
\begin{align*}
f_{1}(x) & =\frac{1}{r}-\frac{\exp(\nu x)}{\beta}\:,\\
f_{2}(x) & =\frac{\rho\mu}{r^{2}}+\frac{\pi(x_{c})}{r}+\frac{\rho}{r}(x-x_{c})\:,
\end{align*}
so that $\mathcal{A}f_{1}(x)+\pi(x)=0$ for $x>x_{c}$ and $\mathcal{A}f_{2}(x)+\pi(x)=0$
for $x<x_{c}$. Then
\[
(R\pi)(x)=\begin{cases}
f_{1}(x)+c_{1}\phi(x) & \text{for }\:x\ge x_{c}\\
f_{2}(x)+c_{2}\psi(x) & \text{for }\:x<x_{c}
\end{cases}\,,
\]
where $c_{1}$ and $c_{2}$ are chosen so that $(R\pi)(x)$ is continuous
and differentiable at $x=x_{c}$. Note that this modification of $\pi(\cdot)$
does not alter the single decision maker solution of Section \ref{subsec:Single}
if $x_{c}$ is sufficiently low (lower than $\theta$ and $x^{*}$)
and if $\rho>rk$ so that Assumption \ref{assump:qi} is satisfied. 

From (\ref{eq:u-x-reg}), we have
\begin{equation}
u(x)=\begin{cases}
0 & \text{for }\:x\ge\theta\\
-\frac{\mathcal{A}V(x)+\pi(x)}{k}=-\mu+\frac{r}{k}V(\theta)+r(x-\theta)-\frac{1}{k}\pi(x) & \text{for}\:x<\theta
\end{cases}\:.\label{eq:u-x}
\end{equation}
Note that $\vert u(x)\vert<\delta(1+\vert x\vert)$ for some constant
$\delta$, so the unique strong solution to the SDE for $Z$ exists. 

\begin{figure}
\begin{centering}
\includegraphics[scale=0.3]{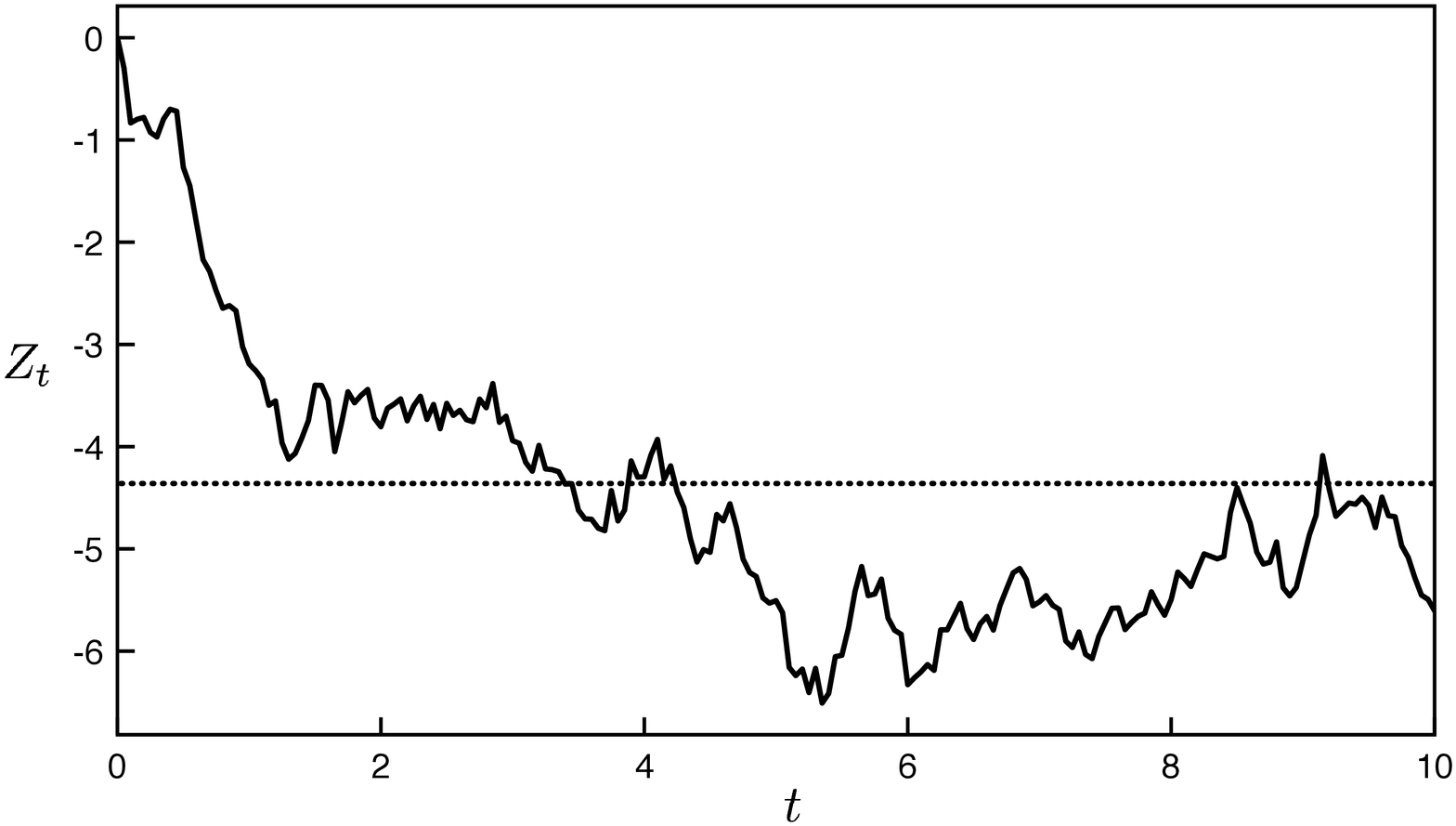}
\par\end{centering}
\begin{centering}
\includegraphics[scale=0.3]{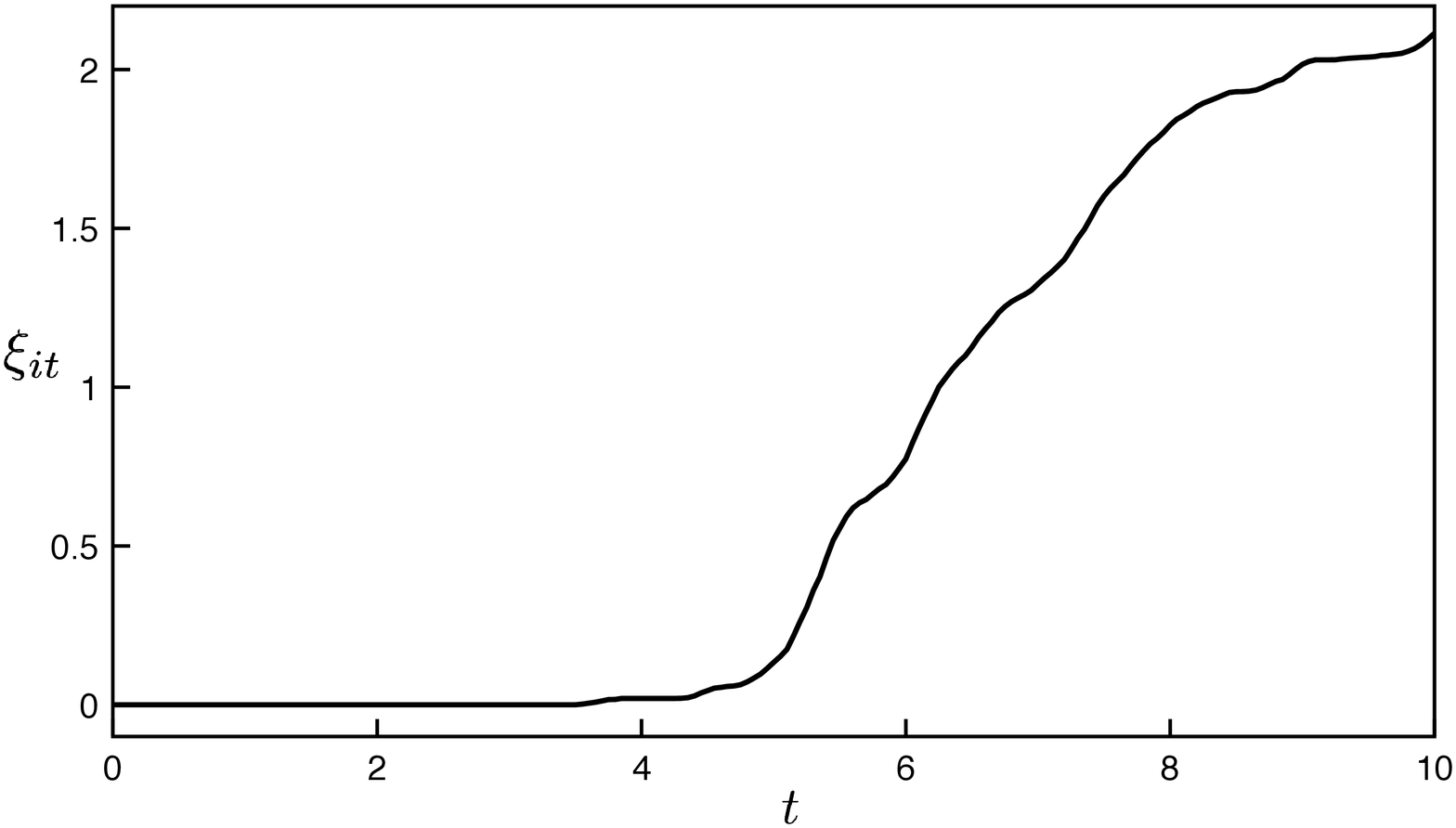}
\par\end{centering}
\begin{centering}
\includegraphics[scale=0.3]{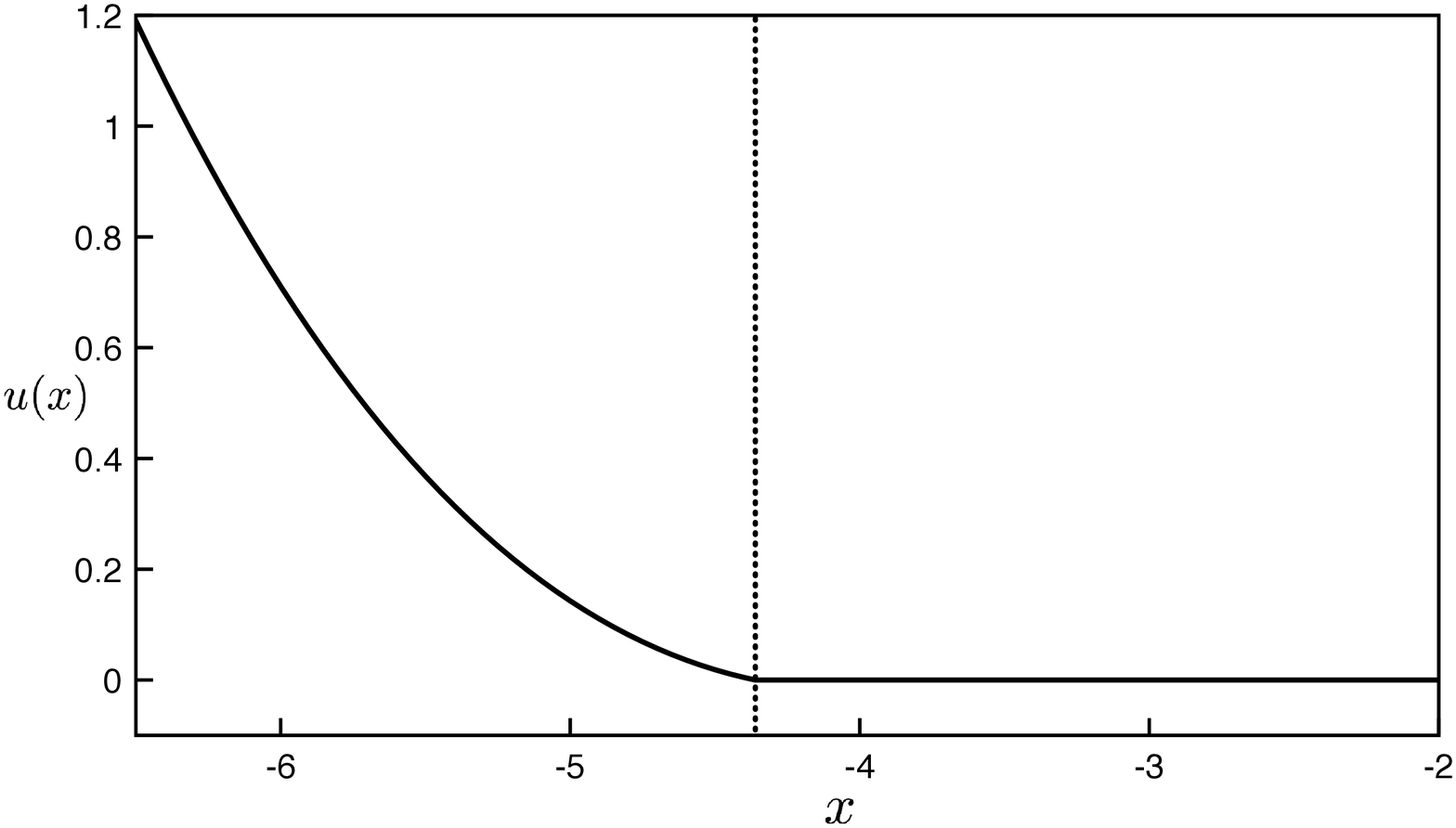}
\par\end{centering}
\caption{Symmetric regular control equilibrium. The dotted lines indicate the
threshold $\theta$. \label{fig:Symm-Eq}}
\end{figure}

Figure \ref{fig:Symm-Eq} illustrates a simulated sample path of $Z$
and $\xi_{i}$ and the rate of contribution $u(\cdot)$ as a function
of $Z$ where $\rho=2$ and $x_{c}=-10$ and the other model parameters
are set as in Figure \ref{fig:Single-decision-maker} including $\theta=-4.3605$.
Note that $Z$ freely fluctuates below $\theta$ in the regular control
equilibrium. In contrast, $Z$ in the single decision maker's solution
shown in Figure \ref{fig:Single-decision-maker} never falls below
$\theta$ because it is subject to singular control at $\theta$.
The sample path of $\xi_{it}=\int_{0}^{t}u(Z_{s})ds$ is smooth and
differentiable with respect to time in the regular control equilibrium,
in contrast to the sample path of $\xi$ in Figure \ref{fig:Single-decision-maker}.
The rate $u(\cdot)$ gradually grows as $Z$ decreases; this is because
the players have stronger incentive to control $Z$ for lower values
of $Z$.

\emph{Remark}: The MPE obtained in Theorem \ref{thm:symm-eq} exists
as long as the solution shown in Lemma \ref{lemma:OptSol} exists.
Hence, Assumptions \ref{assump:mu-sigma}\textendash \ref{assump:qi}
do not need to hold as long as the single decision maker control problem
yields a single-threshold singular control solution.

\section{Impact of Stochasticity and Asymmetry\label{sec:Impact-of-Stoch}}

In this section, we examine the implications of stochasticity and
asymmetry on the regular control equilibrium that we obtain in Section
\ref{subsec:Symmetric-Eq}. This is an important inquiry because most
realistic situations often possess stochastic state variables and
heterogeneity among the players. We first examine the case of deterministic
$Z$ in Section \ref{subsec:Deterministic-Game} and contrast its
equilibrium to that of the stochastic game. We find that the contrast
between a single decision maker solution and an equilibrium is starker
in the stochastic case. Then we examine an asymmetric game in Section
\ref{subsec:Asymmetric-Game} and find that a regular control MPE
does not exist in an asymmetric game. Instead, asymmetry leads to
asymmetric equilibria with singular control strategies.

\subsection{Deterministic Game\label{subsec:Deterministic-Game}}

To examine the impact of stochasticity, we will simply discuss the
deterministic case of Example 2 and contrast it to the stochastic
case.

\emph{Example 3}: We revisit the model of Example 2 and set $\sigma=0$
and $x_{c}=-\infty$. Then the characteristic operator
\[
\mathcal{A}=\mu\partial_{z}-r
\]
is a first-order differential operator, and there is only one fundamental
solution $\phi(x)=\exp(rx/\mu)$ that satisfies $\mathcal{A}\phi(x)=0$.
Furthermore, 
\[
(R\pi)(x)=\frac{1}{r}-\frac{\exp(\nu x)}{r-\mu\nu}\:,
\]
where we assume $\mu\nu<r$. The critical threshold of control is
given by
\[
\theta=\frac{1}{\nu}\ln\frac{kr}{\vert\nu\vert}\;.
\]
 As in Example 2, $u(\cdot)$ is given by (\ref{eq:u-x}). It is straightforward
to verify that $\lim_{x\uparrow\theta}u(x)=0$, $\lim_{x\uparrow\theta}u'(x)=0$
and that $u''(x)>0$ for all $x<\theta$, which implies that $u'(x)<0$
for all $x<\theta$. Furthermore, there exists $\eta<\theta$ at which
$\mu+2u(\eta)=0$ so that $dZ_{t}/dt<0$ for $Z_{t}>\eta$ and $dZ_{t}/dt>0$
for $Z_{t}<\eta$. 

\begin{figure}
\begin{centering}
\includegraphics[scale=0.4]{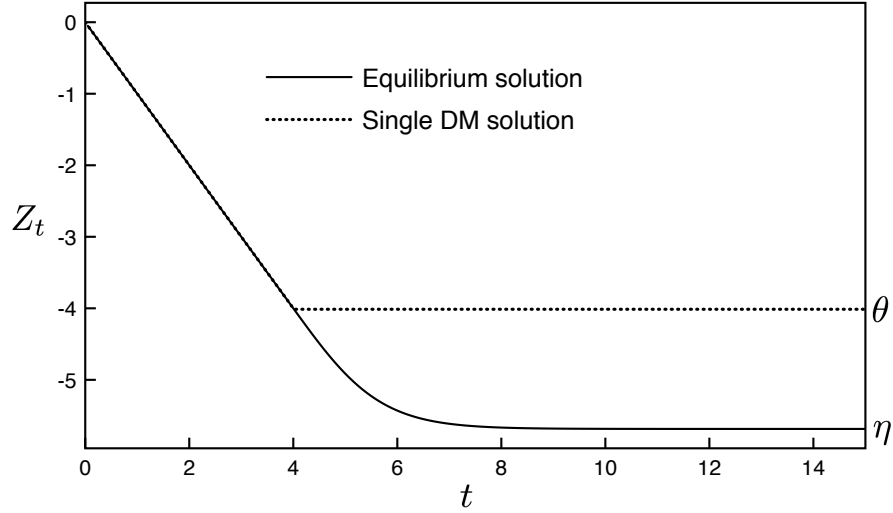}
\par\end{centering}
\caption{Deterministic solutions. \label{fig:Deterministic}}
\end{figure}

Lastly, we can show that $Z_{t}$ asymptotically approaches $\eta$
if $Z_{0}\neq\eta$. Define $y_{t}=Z_{t}-\eta$. For $y_{t}$ sufficiently
close to 0, we have 
\begin{align*}
\frac{dy_{t}}{dt} & =\mu+2u(y_{t}+\eta)=\mu+2u(\eta)+2u'(\eta)y_{t}+O(y_{t}^{2})\\
 & =2u'(\eta)y_{t}+O(y_{t}^{2})\:,
\end{align*}
where we used the fact that $\mu+2u(\eta)=0$ in the last equality.
Recall that $u'(\eta)<0$; thus, $\vert y_{t}\vert=\vert y_{0}\vert\exp[2u'(\eta)t+O(y_{t})]$
for large $t$, which implies $\lim_{t\rightarrow\infty}y_{t}=0$.
Hence, $\eta$ is the steady state of $Z$. Figure \ref{fig:Deterministic}
illustrates $Z_{t}$ as a function of $t$ in the equilibrium where
$\theta=-4.0132$ and $\eta=-5.68646$. 

In the deterministic case of the single decision maker problem, $Z_{t}$
follows the horizontal line $z=\theta$ as soon as $Z_{t}$ hits $\theta$.
Thus, $\theta$ is the steady state of the single decision maker solution.
Hence, in the deterministic model, the behavior of $Z_{t}$ exhibits
very little qualitative difference between the single decision maker
solution and the equilibrium of the game. In particular, an outside
observer will not be able to tell the difference between the two solutions
except that the steady state values differ. The difference in the
steady state value of $Z$ can be simply attributed to the free rider
effect: the players are less willing to contribute to the common good,
so the steady state is lower.

In contrast, in the presence of stochasticity ($\sigma>0$), the behavior
of $Z_{t}$ is markedly different as illustrated in Figures \ref{fig:Single-decision-maker}
and \ref{fig:Symm-Eq}: in the single decision maker's case, $Z_{t}$
is reflected off of the threshold $\theta$, whereas in the equilibrium
of the game, $Z_{t}$ can assume any value although it tends to fluctuate
around $\theta$. We conclude that stochasticity induces qualitatively
different behaviors between the two solutions and hence renders the
regular control equilibrium observable to an outsider. 

\subsection{Asymmetric Game \label{subsec:Asymmetric-Game}}

Next, we examine the impact of asymmetry between the two players.
We first show that a regular control MPE is absent in an asymmetric
game in Section \ref{subsubsec:Absence-of-a}, and then we construct
the simplest class of asymmetric equilibria in Section \ref{subsubsec:Asymmetric-Eq}
and demonstrate that they exhibit the key characteristics of a singular
control solution.

\subsubsection{Absence of a Regular Control MPE \label{subsubsec:Absence-of-a}}

Suppose that $\theta_{i}\neq\theta_{j}$ because of asymmetry ($k_{i}\neq k_{j}$
and/or $\pi_{i}\neq\pi_{j}$), where $\theta_{i}$ is the unique solution
to the equation
\[
\frac{k_{i}-(R\pi_{i})'(\theta_{i})}{\phi'(\theta_{i})}=-\frac{(R\pi_{i})''(\theta_{i})}{\phi''(\theta_{i})}\:\:.
\]
For analytical tractability, we make the following additional assumption:

\begin{assumption}\label{assum:pi-p-rk} $\mu'(x)+\pi'(x)-rk\neq0$
almost everywhere $x\in I$. \end{assumption} Assumption \ref{assum:pi-p-rk}
ensures that $\mu(x)+\pi(x)-rk$ is never constant within any given
non-empty interval. It gives a non-trivial structure to the Hamilton-Jacobi-Bellman
(HJB) equation of the payoff function. 

The following theorem establishes that there is no payoff function
$V_{i}(\cdot;\xi)\in C^{2}(I)$ associated with a regular control
MPE.

\begin{theorem} \label{thm:asymmetric} If $\theta_{i}\neq\theta_{j}$,
then there is no regular control MPE $\xi=(\xi_{i},\xi_{j})\in\Sigma_{i}\times\Sigma_{j}$
such that $V_{i}(\cdot;\xi)\in C^{2}(I)$ for $i=1,2$. \end{theorem} 

The implication of Theorem \ref{thm:asymmetric} is that an MPE of
an asymmetric game must involve singular control. Hence, the natural
next step is to explore the form of such a singular control MPE in
an asymmetric game, which is the goal of Section \ref{subsubsec:Asymmetric-Eq}.

\emph{Remark 1}: Strictly speaking, Theorem \ref{thm:asymmetric}
does not necessarily preclude the possibility of an equilibrium with
payoffs $V_{i}(\cdot;\xi)$ that do not belong to $C^{2}(I)$ such
as a general \emph{viscosity solution}. In this paper, we limit ourselves
to equilibria that produce \emph{classic solutions} only and defer
the possibility of an equilibrium with non-$C^{2}(I)$ viscosity solutions
to future endeavors. We also remark that we do not attempt to exclude
the possibility of an equilibrium $\xi$ that does not belong to $\Sigma_{i}\times\Sigma_{j}$.
This possibility is beyond the scope of the paper because of the issue
of the existence of the unique solution to the SDE for $Z_{t}$; this
is a common restriction in stochastic control theory \citep{Oksendal2005}. 

\emph{Remark 2}: Even in an asymmetric game, it is possible to have
some combinations of $k_{i}$ and $\pi_{i}(\cdot)$ such that $\theta_{1}=\theta_{2}$.
In this case, a regular control MPE is possible because Theorem \ref{thm:symm-eq}
is applicable. 

As a corollary, we can also exclude the possibility of an MPE in which
there is a common regular control region $(\alpha,\beta)$ for some
$\alpha>a$ and a no-control region $[\beta,b)$. (We remark that
the proof of Theorem \ref{thm:asymmetric} establishes that the regular
control regions of the two players must coincide.) We remain agnostic
about what happens in the region $(a,\alpha]$, however, except that
$(a,\alpha]$ should contain singular control regions $D_{1}$ and
$D_{2}$ if they exist. Again, we focus on an MPE with classical solutions,
i.e., $V_{i}(\cdot;\xi)\in C^{2}(I\backslash D_{j})$. 

\begin{corollary} \label{cor:preclude} There is no MPE $\xi\in\Sigma_{i}\times\Sigma_{j}$
such that $V_{i}(\cdot;\xi)\in C^{2}(I\backslash D_{j})$ with a common
regular control region $(\alpha,\beta)$ for some $\alpha>a$ and
a no-control region $[\beta,b)$. \end{corollary} Its proof essentially
follows that of Theorem \ref{thm:asymmetric}, so it is omitted. 

\subsubsection{Non-regular Control Asymmetric MPE \label{subsubsec:Asymmetric-Eq}}

By virtue of Theorem \ref{thm:asymmetric}, an asymmetric game allows
no regular control MPE, so any MPE must involve \emph{some} singular
control by at least one player. Furthermore, Corollary \ref{cor:preclude}
implies that the only possible equilibria are the ones with a singular
control region $[\alpha,\beta]$ of one player and a no-control region
$(\beta,b)$. Our goal is to present the simplest class of such equilibria
and compare its characteristics to those of the regular control MPE.
Although our main focus is on strictly asymmetric games, we will keep
our discourse general and implicitly include the case of a symmetric
game.

As a candidate for an asymmetric MPE, we consider a strategy profile
$\xi$ in which player 1 exerts singular control in the interval $[\theta',\theta_{1}]$
and regular control in $(a,\theta')$ for some $\theta'<\theta_{1}$
while player 2 exerts regular control in $(a,\theta')$. Furthermore,
we define the regular control rate functions
\begin{align}
u_{1}(x) & =\begin{cases}
-\frac{1}{k_{2}}[\mathcal{A}U_{2}(x)+\pi_{2}(x)] & \text{for}\:x<\theta'\:,\\
0 & \text{otherwise}
\end{cases}\:,\label{eq:u1-asymm}\\
u_{2}(x) & =\begin{cases}
-\frac{1}{k_{1}}[\mathcal{A}U_{1}(x)+\pi_{1}(x)] & \text{for}\:x<\theta'\:,\\
0 & \text{otherwise}
\end{cases}\:,\label{eq:u2-asymm}
\end{align}
where $U_{1}(\cdot)$, $U_{2}(\cdot)$ are given by

\begin{align*}
U_{1}(x) & =V_{1}(x)\;,\\
U_{2}(x) & =\begin{cases}
B\phi(x)+(R\pi_{2})(x) & x\ge\theta_{1}\\
U_{2}(\theta_{1}) & x\in(\theta',\theta_{1})\\
U_{2}(\theta_{1})+(x-\theta')k_{2} & x\le\theta'
\end{cases}\:,
\end{align*}
Here $V_{1}(\cdot)$ denotes the solution to the single decision maker
problem given by (\ref{eq:Vx}) where $k$, $\theta$ and $\pi$ are
respectively replaced by $k_{1}$, $\theta_{1}$ and $\pi_{1}$, and
\begin{align}
B & =-\frac{(R\pi_{2})'(\theta_{1})}{\phi'(\theta_{1})}\;,\label{eq:B-theta}\\
U_{2}(\theta_{1}) & =B\phi(\theta_{1})+(R\pi_{2})(\theta_{1})\:.\nonumber 
\end{align}

Note that the functions $U_{1}(\cdot)$, $U_{2}(\cdot)$ are actually
the payoff functions associated with the proposed strategy profile
$\xi$. In the region $(\theta_{1},b)$, both $U_{1}(\cdot)$ and
$U_{2}(\cdot)$ assume the form of continuation region without control.
In player 1's singular control region $[\theta',\theta_{1}]$, we
have $U_{1}'(x)=k_{1}$ because player 1 expends the cost of singular
control while $U_{2}'(x)=0$ because player 2 does not expend any
cost in this region. In the common regular control region $(a,\theta')$,
both players expend cost in such a way that $U_{i}'(x)=k_{i}$. 

By Theorem \ref{thm:asymmetric}, to confirm that the strategy profile
$\xi$ is an MPE, we only need to verify the following sufficient
conditions:
\begin{align}
U_{2}'(x) & \le k_{2}\:\forall x\in I\backslash\{\theta_{1},\theta'\},\label{eq:U2'}\\
\mathcal{A}U_{2}(x)+\pi_{2}(x) & \le0\;\forall x\in I\backslash\{\theta_{1},\theta'\},\label{eq:AU2}\\
u_{1}(\cdot)\in\mathcal{U}\,,\: & u_{2}(\cdot)\in\mathcal{U}\quad.\label{eq:U}
\end{align}
 \begin{proposition} \label{prop:asymm-eq} Suppose that Assumptions
\ref{assump:mu-sigma}\textendash \ref{assump:qi} hold and (\ref{eq:U2'})-(\ref{eq:U})
are satisfied. Then the strategy profile $\xi$ is an MPE with the
payoff functions given by $V_{1}(x;\xi)=U_{1}(x)$ and $V_{2}(x;\xi)=U_{2}(x)$.
\end{proposition} 

\begin{figure}
\begin{centering}
\includegraphics[scale=0.4]{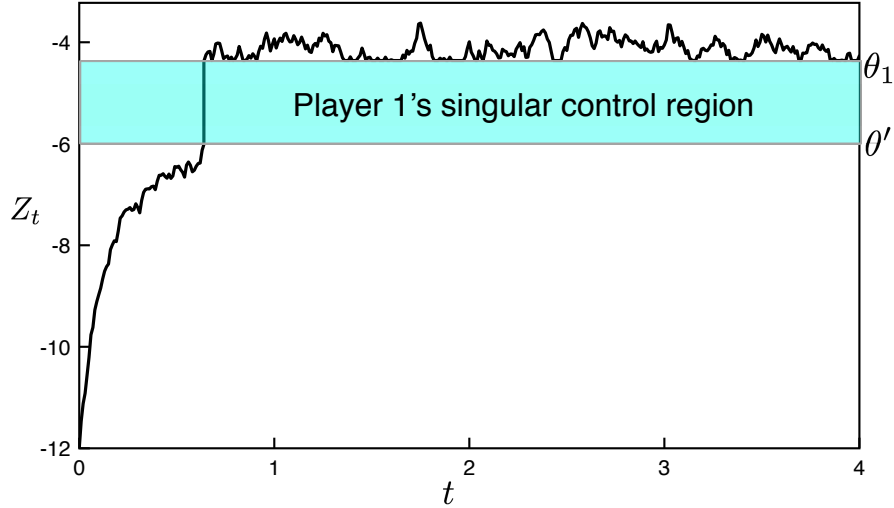}
\par\end{centering}
\caption{Sample path of $Z_{t}$ in an asymmetric equilibrium. \label{fig:asymmetric}}
\end{figure}

\emph{Example 4}: Recall Example 2 and set $\mu=-1$, $\sigma=r=k_{1}=k_{2}=1$,
$\rho=2$, and $x_{c}=-10$ just as in Figure \ref{fig:Symm-Eq}.
Then $\theta=-4.3605$ is the boundary of the control region. Here
we can construct an asymmetric equilibrium as in Proposition \ref{prop:asymm-eq}
because we can verify that any choice of $\theta'<\theta_{1}=\theta$
satisfies conditions (\ref{eq:U2'}) and (\ref{eq:AU2}). Hence, there
is a continuum of asymmetric equilibria parameterized by $\theta'\in(-\infty,\theta_{1})$.
Figure \ref{fig:asymmetric} illustrates a numerical example of the
case $\theta'=-6$. Note that even if we set $k_{1}<k_{2}$, the qualitative
features of this asymmetric equilibrium continue to hold.

In particular, Figure \ref{fig:asymmetric} illustrates the case when
$Z_{0}<\theta'$. Before $Z_{t}$ hits $\theta'$ for the first time,
$Z_{t}$ is subject to regular control of both players. Upon reaching
$\theta'$, $Z_{t}$ is subject to singular control by player 1 and
is boosted up to $\theta_{1}$. Once $Z_{t}$ enters the region $[\theta_{1},\infty)$,
the threshold $\theta_{1}$ takes the role of a reflecting boundary
for $Z_{t}$ because of player 1's singular control strategy. If $Z_{t}\in[\theta_{1},\infty)$,
the equilibrium reduces to the single decision maker solution.

As illustrated by the example, asymmetric non-regular equilibria exist
even for the symmetric game where $k_{1}=k_{2}$ and $\pi_{1}=\pi_{2}$.
However, if the two players are identical to each other, the players
are likely to be drawn to the symmetric equilibrium. For instance,
suppose that the equilibrium shown in Figure \ref{fig:asymmetric}
is the outcome of the symmetric game. Then eventually player 1 ends
up being the only one contributing to the common good every time $Z_{t}$
hits $\theta$, so he will feel that the current equilibrium is unfair
to him. Consequently, he will likely attempt to switch to a more equitable
equilibrium. Thus, the symmetric regular MPE is the likely focal point
equilibrium \citep{Fudenberg1991}. In contrast, in an asymmetric
game in which the two players have unequal thresholds $\theta_{i}\neq\theta_{j}$,
the only possible equilibrium is an asymmetric non-regular control
one by virtue of Theorem \ref{thm:asymmetric}. 

One interesting question regards which player exerts singular control.
In the numerical example above, we can fix $k_{1}=1$ and vary the
value of $k_{2}$ and see how that affects the equilibrium. It can
be numerically verified that the MPE of Proposition \ref{prop:asymm-eq}
exists as long as $k_{2}\ge0.5383$. If $k_{2}$ is less than the
critical value 0.5383, then (\ref{eq:U2'}) is violated, so the MPE
is not possible. Intuitively, if $k_{2}$ is sufficiently low, then
player 2 has strong incentive to exert singular control, and knowing
this, player 1 would never exert singular control. In this case, the
only equilibrium is the one in which player 2 exerts singular control.
Thus, sufficiently high asymmetry induces the more efficient player
to exert singular control in an equilibrium. If the asymmetry is modest,
then either player can be the one who exerts singular control.

In summary, the MPE obtained by Proposition \ref{prop:asymm-eq} is
the simplest form of asymmetric equilibria. In these equilibria, $Z_{t}$
eventually ends up in the region $[\theta_{1},b)$ and subject to
the reflecting boundary at $\theta_{1}$. Thus, in the long run, $Z_{t}$
is subject to the singular control policy of player 1 just as in the
single decision maker's case, and it is not plagued by the gradualism
of the regular control MPE of Section \ref{subsec:Symmetric-Eq}.
Therefore, we conclude that asymmetry reduces inefficiency. 

\section{Discussions \label{sec:Discussions}}

\subsection{Dimensionality of the State Variable}

In contrast to the singular control equilibria obtained in, for example,
\citet{Ghemawat1990}, \citet{Steg2012}, \citet{Battaglini2014},
\citet{Ferrari2017}, and Appendix \ref{sec:R=000026D-Game}, a regular
control equilibrium arises in our model. We speculate that we obtain
a contrasting result because of the difference in the dimensionality
of the state variable. In the model that we study, the state variable
$Z_{t}$ is one-dimensional; the control variables $\xi_{it}$ and
$\xi_{jt}$ only add to $Z_{t}$, so they are not independent state
variables that stand alone. In contrast, in the examples from the
literature as well as the R\&D spillover game analyzed in Appendix
\ref{sec:R=000026D-Game}, the state variables are multidimensional
because the players' control variables are decoupled from the state
variable. For instance, in the R\&D game of Appendix \ref{sec:R=000026D-Game},
the state variable is two-dimensional: $(\lambda_{1t},\lambda_{2t})$,
where $\lambda_{it}$ is the current level of R\&D effort of firm
$i$. Consequently, the possibility of a subgame with $\lambda_{1t}\neq\lambda_{2t}$
is allowed, so the state variable is allowed to be asymmetric between
the two players. However, as shown by Section \ref{subsubsec:Asymmetric-Eq},
the regular control equilibrium of Theorem \ref{thm:symm-eq} hinges
on the symmetry of the state variable between the two players so that
they share the common regular control region (see the proof of Theorem
\ref{thm:symm-eq}). Thus, we anticipate that the emergence of a regular
control equilibrium is driven by the single-dimensionality of the
state variable.

\subsection{N-Player Game \label{subsec:N-Player-Game}}

We can straightforwardly generalize Theorem \ref{thm:symm-eq} to
an $N$-player game for $N>2$ by constructing a symmetric regular
control equilibrium. We first define $u^{N}(\cdot)$ as follows:
\[
u^{N}(x)=\begin{cases}
-\frac{1}{k(N-1)}\left(\mathcal{A}V(x)+\pi(x)\right) & \text{for}\;x<\theta\\
0 & \text{for}\;x\ge\theta
\end{cases}\:,
\]
where $V(\cdot)$ is given by (\ref{eq:Vx}). Then it can be verified
that a strategy profile $\xi$ in which each player $i$ exerts a
regular control of $d\xi_{it}=u^{N}(Z_{t})dt$ constitutes a MPE.
This is because the HJB equation (condition (iii) of Theorem \ref{thm:verif})
for each player's payoff function continues to be satisfied when all
$N-1$ opponents exert the regular control of $u^{N}(Z_{t})dt$. 

\subsection{Relation to Mixed Strategy Equilibrium of War of Attrition \label{subsec:Relation-to-Mixed}}

There exists a close analogy between the regular control equilibrium
obtained in Theorem \ref{thm:symm-eq} and the mixed strategy equilibrium
of a war of attrition. In the mixed strategy equilibrium of the canonical
war of attrition, each player has control over the hazard rate of
exit, which is analogous to the rate of regular control in our model.
Similarly, a pure strategy concession (a deterministic concession)
in a war of attrition is analogous to the singular control strategy
in our model. 

The analogy between the two equilibrium solutions goes even further
with the impact of asymmetry. In a stochastic extension of a war of
attrition game, the mixed strategy MPEs disappear when the players's
reward from concession is asymmetric \citep{Georgiadis2019}. This
is analogous to our result that a completely regular control MPE disappears
when the players are asymmetric. 

\section{Conclusions \label{sec:Conclusions}}

We examine a stochastic game of variable contribution as a generalization
of a war of attrition. In particular, we analyze a stochastic game-theoretic
extension of the Nerlove-Arrow model, which possesses a novel MPE
characterized by regular control. This finding is in contrast to the
singular control equilibria possessed by variable concession games
with multidimensional state variables. In the examples of singular
control equilibria obtained in the literature, the free rider effect
manifests in the value of the threshold of the control region, but
the action of concession is immediate and not plagued by delay or
gradualism. In contrast, the variable contribution game analyzed in
Section \ref{sec:VariableContribution} possesses a regular control
equilibrium in which the free rider effect manifests in the gradualism
of the players' actions. 

We find that it is important to understand the effect of stochasticity
on the game. The state variable $Z$ exhibits qualitatively different
behavior under a regular control MPE from that of a single decision
maker solution. However, the difference almost disappears if the state
variable is deterministic. We conclude that stochasticity renders
the gradual MPE observable to an outsider. 

We also examine the impact of asymmetry between the players and find
that the regular control MPE is not possible under asymmetry. The
implication of this finding is that the problem of inefficiency arising
from the gradual regular control MPE is mitigated by asymmetry between
the players. From a social planner's perspective, this result suggests
that heterogeneity between agents should be cultivated or encouraged
when there is a free rider problem with the agents' contributions
to the common good. 

The results and their implications of this paper warrant some related
future research endeavors. First, it will be interesting to study
multidimensional variable concession problems and see if they have
gradual regular control MPE even though we speculate that they do
not. Second, it will be fruitful to examine an extension of our model
in which the players have private types and asymmetric information
regarding the cost of contribution. In this case, there is inherent
asymmetry between any two players, so the regular control equilibrium
may exist only under very stringent conditions. Lastly, just as \citet{Wang2009}
finds empirical evidence of a delay in action in a war of attrition,
it might be possible to find empirical evidence of gradualism in a
regular control equilibrium for a contribution game with a free rider
problem such as in generic advertising or investment in public goods.

{\small{}\bibliographystyle{INFORMS2011}
\bibliography{GradualWA}
}{\small\par}

\newpage{}

\appendix

\section{Mathematical Proofs \label{sec:Mathematical-Proofs}}

\textbf{Proof of Lemma \ref{lemma:OptSol}}: Because the first and
second derivatives of $V(\cdot)$ are continuous at $\theta$, the
coefficient $A$ and the threshold $\theta$ must satisfy $k=(R\pi)'(\theta)+A\phi'(\theta)$
and $0=(R\pi)''(\theta)+A\phi''(\theta)$. The simultaneous equations
are solved if we can find $\theta$ that satisfies 
\begin{equation}
\frac{k-(R\pi)'(\theta)}{\phi'(\theta)}=-\frac{(R\pi)''(\theta)}{\phi''(\theta)}\:.\label{eq:th-nece}
\end{equation}
Below we show that there is a unique value of $\theta$ that satisfies
this condition and that the resulting solution $V(\cdot)$ satisfies
the optimality conditions for singular stochastic control. 

As a first step, we show that the function $F(x):=(k-(R\pi)'(x))/\phi'(x)$
achieves a unique global maximum at $x=\theta\in I$. If this holds,
then it is straightforward to verify that (\ref{eq:th-nece}) is satisfied.
Note that $(Rq)'(x)=(R\pi)'(x)-k$ from the definition of $q$ in
(\ref{eq:q1-x-1}) without the player index $i$, so $F(x)=-(Rq)'(x)/\phi'(x)$.
From the theory of diffusive processes \citep{Borodin1996,Alvarez2008},
it is well-known that 
\[
(Rq)(x)=B^{-1}\left[\phi(x)\int_{a}^{x}\psi(y)q(y)m'(y)dt+\psi(x)\int_{x}^{b}\phi(y)q(y)m'(y)dt\right]\:,
\]
where 
\begin{align*}
B & =\frac{1}{S'(x)}[\psi'(x)\phi(x)-\phi'(x)\psi(x)]\;,\\
m'(x) & =\frac{2}{\sigma^{2}(x)S'(x)}\quad,\\
S'(x) & =\exp\left(-\int_{0}^{x}\frac{2\mu(y)}{\sigma^{2}(y)}dy\right)\;.
\end{align*}
 Thus, we have
\begin{align}
\frac{dF(x)}{dx} & =-\frac{2S'(x)}{\sigma^{2}(x)[\phi'(x)]^{2}}L(x)\:,\label{eq:dF-dx}\\
\text{where}\;L(x) & =-\frac{q(x)\phi'(x)}{S'(x)}-r\int_{x}^{b}\phi(y)q(y)m'(y)dy\:.-x\label{eq:L-x}
\end{align}
 By the definitions of $m'(\cdot)$ and $S'(\cdot)$, we can derive
the equality $d(\phi'(x)/S'(x))/dx=r\phi(x)m'(x)$ by some algebra,
from which we obtain 
\[
\frac{dL(x)}{dx}=-q'(x)\frac{\phi'(x)}{S'(x)}\:.
\]
 By Assumption \ref{assump:qi} (ii), $q'(x)>0$ for $x\in(a,x^{*})$
and $q'(x)<0$ for $(x^{*},b)$, we have $dL(x)/dx>0$ for $x\in(a,x^{*})$
and $dL(x)/dx<0$ for $x\in(x^{*},b)$. From $\lim_{x\rightarrow a}\phi'(x)/S'(x)=-\infty$
(p. 19 of \citet{Borodin1996}) and $\lim_{x\rightarrow a}q(x)=-\infty$
from Assumption \ref{assump:qi} (iii), we have for $y$ sufficiently
close to $a$, 
\[
L(y)-\lim_{x\rightarrow a}L(x)=\lim_{x\rightarrow a}\int_{x}^{y}\frac{dL(u)}{du}du=\lim_{x\rightarrow a}\int_{x}^{y}[-q'(u)]\frac{\phi'(u)}{S'(u)}du>\frac{\phi'(y)}{S'(y)}\lim_{x\rightarrow a}\int_{x}^{y}[-q'(u)]du=\infty\:,
\]
 so that $\lim_{x\rightarrow a}L(x)=-\infty$. Furthermore, because
$\lim_{x\rightarrow b}q(x)\phi'(x)/S'(x)=0$ from Assumption \ref{assump:qi}
(iii), we also have $\lim_{x\rightarrow b}L(x)=0$ by the expression
(\ref{eq:L-x}). From the continuity of $L(\cdot)$ and the sign change
of $dL(x)/dx$, it follows that $L(\theta)=0$ for some unique point
$\theta\in(a,x^{*})$. Since $L(\cdot)$ is monotonically increasing
in the interval $(a,x^{*})$, $L(\cdot)$ turns from negative to positive
at $\theta$. Combined with the behavior of $L(x)$ in the limits
$x\rightarrow a,b$, we conclude that $L(x)$ is negative in the interval
$(a,\theta)$ and positive in the interval $(\theta,b)$. From (\ref{eq:dF-dx}),
we also conclude that $F(x)$ attains its global maximum at this unique
point $\theta$. Thus, (\ref{eq:th-nece}) is also satisfied, which
makes $V(\cdot)\in C^{2}(I)$ if $A=F(\theta)$. 

The next step is to prove that $V(\cdot)$ satisfies the sufficient
conditions (i) and (ii) of the optimality. 

(i) First, we show that $\mathcal{A}V(x)+\pi(x)\le0$ and $V'(x)-k\le0$
for all $x\in I$. For $x>\theta$, the form of $V(\cdot)$ guarantees
the condition $\mathcal{A}V(x)+\pi(x)=0$. Furthermore, because $F(x)<F(\theta)=A$
for all $x>\theta$, we have $(k-(R\pi)'(x))/\phi'(x)<A$, from which
we obtain
\[
V'(x)=A\phi'(x)+(R\pi)'(x)<k
\]
 for all $x>\theta$.

For $x<\theta$, we have $V'(x)=k$ by the form of $V(\cdot)$. Also
note that $V''(\theta)=0$, $V'(\theta)=k$ and $\mathcal{A}V(x)+\pi(x)=0$
for $x>\theta$ so that 
\[
\lim_{x\downarrow\theta}\mathcal{A}V(x)+\pi(x)=\mu(\theta)k-rV(\theta)+\pi(\theta)=0\:.
\]
For any $x<\theta$, we have
\begin{align*}
\mathcal{A}V(x)+\pi(x)-[\lim_{y\downarrow\theta}\mathcal{A}V(y)+\pi(y)] & =\mu(x)k-rV(x)+\pi(x)-[\mu(\theta)k-rV(\theta)+\pi(\theta)]\\
 & =q(x)-q(\theta)<0
\end{align*}
where the inequality is from the fact that $q(\cdot)$ increases in
the interval $(a,\theta)$. Thus, $\mathcal{A}V(x)+\pi(x)<0$ for
$x<\theta$.

(ii) We just showed that $\mathcal{A}V(x)+\pi(x)=0$ for $x\ge\theta$
and $V'(x)-k=0$ for $x\le\theta$. \eproof 

\textbf{Proof of Theorem \ref{thm:verif}}: To prove the theorem,
it is sufficient to show that $U(x)=V_{i}(x;\xi_{i}^{*},\xi_{j})\ge V_{i}(x;\xi_{i},\xi_{j})$
for any arbitrary strategy $\xi_{i}$ of player $i$ that satisfies
$d\xi_{it}=u_{i}(t,Z_{t})+d\xi_{it}^{c}+\Delta\xi_{it}$. First, it
is straightforward to verify that if $\xi_{i}^{*}$ is the best response,
player $i$ should not expend any cost to control $Z$ within $D_{j}$
because player $j$ is already doing so in that region. Thus, one
necessary condition for the payoff function of the best response is
$U'(x)=0$ in the interior of $D_{j}$. 

We consider an arbitrary strategy $\xi_{i}$ of player $i$ that satisfies
$d\xi_{it}=u_{i}(t,Z_{t})dt+d\xi_{it}^{c}+\Delta\xi_{it}$ for some
arbitrary $u_{i}(\cdot)\in\mathcal{U}$. Let $Z$ be the state process
dictated by the given strategy profile $(\xi_{i},\xi_{j})$. 

By conditions (i) and (iii), $U'(\cdot)$ is a bounded function. Furthermore,
due to Assumption \ref{assump:mu-sigma} (i), $e^{-rs}\sigma(Z_{s})\tilde{U}'(Z_{s})$
is locally bounded. Hence, the process $M_{t}:=\int_{0}^{t}e^{-rs}\sigma(Z_{s})\tilde{U}'(Z_{s})dW_{s}$
is a continuous local martingale. By the definition of a continuous
local martingale \citep{Karatzas1998}, there exists a non-decreasing
sequence $\{\tau_{n}\}_{n=1}^{\infty}$ of stopping times of $\{\mathcal{F}_{t}\}_{t>0}$
such that $\lim_{n\rightarrow\infty}\tau_{n}=\infty$ a.s. and $\{M_{\tau_{n}\wedge t}:t\in[0,\infty)\}$
is a martingale for each $n$. We first consider any $x\in C_{j}$
as the initial point of $Z$. By the generalized It\^{o}'s formula,
we have
\begin{align*}
e^{-r\tau_{n}\wedge t}\tilde{U}(Z_{\tau_{n}\wedge t})= & \tilde{U}(x)+\int_{0}^{\tau_{n}\wedge t}e^{-rs}\mathcal{A}\tilde{U}(Z_{s})ds\\
 & +\int_{0}^{\tau_{n}\wedge t}e^{-rs}\tilde{U}'(Z_{s})[u_{i}(t,Z_{s})ds+u_{j}(Z_{s})ds+\sum_{l=1}^{2}d\xi_{ls}^{c}]\\
 & +\sum_{0\le s\le\tau_{n}\wedge t}e^{-rs}[\tilde{U}(Z_{s})-\tilde{U}(Z_{s^{-}})]+M_{\tau_{n}\wedge t}\:.
\end{align*}
Taking the expectation of both sides and rearranging terms, we obtain
\begin{align}
\tilde{U}(x)= & \mathbb{E}^{x}[e^{-r\tau_{n}\wedge t}\tilde{U}(Z_{\tau_{n}\wedge t})]-\mathbb{E}^{x}\left[\int_{0}^{\tau_{n}\wedge t}e^{-rs}(\mathcal{A}\tilde{U}(Z_{s})+\tilde{U}'(Z_{s})u_{j}(Z_{s})+\tilde{U}'(Z_{s})u_{i}(s,Z_{s}))ds\right]\nonumber \\
 & -\mathbb{E}^{x}\left[\int_{0}^{\tau_{n}\wedge t}e^{-rs}\tilde{U}'(Z_{s})\sum_{l=1}^{2}d\xi_{ls}^{c}+\sum_{0\le s\le\tau_{n}\wedge t}e^{-rs}[\tilde{U}(Z_{s})-\tilde{U}(Z_{s^{-}})]\right]\:.\label{eq:U-x}
\end{align}
Here we use the fact that $\mathbb{E}^{x}[M_{\tau_{n}\wedge t}]=0$
because $M_{\tau_{n}\wedge t}$ is a martingale. 

Next, we re-express (\ref{eq:U-x}) as an inequality involving $U(\cdot)$
and $\xi_{i}$ alone. Recall that $\tilde{U}(x)=U(x)$ whenever $x\in C_{j}$,
so $\tilde{U}(Z_{t})=U(Z_{t})$ whenever $Z_{t}$ is continuous in
time. Thus, if $Z_{s}\in C_{j}$, 
\[
\mathcal{A}\tilde{U}(Z_{s})+\tilde{U}'(Z_{s})(u_{i}(t,Z_{s})+u_{j}(Z_{s}))=\mathcal{A}U(Z_{s})+U'(Z_{s})(u_{i}(t,Z_{s})+u_{j}(Z_{s}))\,.
\]
We also note that $\tilde{U}'(Z_{s})d\xi_{js}^{c}=0$ for all $s$
from condition (i) since $d\xi_{js}^{c}=0$ if $Z_{s}\in D_{j}$ and
$\tilde{U}'(Z_{s})=U'(Z_{s})$ if $Z_{s}\in C_{j}$. Also note that
the process $Z_{t}$ spends zero time within $D_{j}$, so $\tilde{U}(Z_{s})-\tilde{U}(Z_{s^{-}})=U(Z_{s})-U(Z_{s^{-}})$
for all $s>0$. Then we can re-express (\ref{eq:U-x}) as the following
equality: 
\begin{align}
U(x)= & \mathbb{E}^{x}[e^{-r\tau_{n}\wedge t}U(Z_{\tau_{n}\wedge t})]-\mathbb{E}^{x}\left[\int_{0}^{\tau_{n}\wedge t}e^{-rs}(\mathcal{A}U(Z_{s})+U'(Z_{s})u_{j}(Z_{s})+U'(Z_{s})u_{i}(s,Z_{s}))ds\right]\nonumber \\
 & -\mathbb{E}^{x}\left[\int_{0}^{\tau_{n}\wedge t}e^{-rs}U'(Z_{s})\sum_{l=1}^{2}d\xi_{ls}^{c}+\sum_{0\le s\le\tau_{n}\wedge t}e^{-rs}[U(Z_{s})-U(Z_{s^{-}})]\right]\:.\label{eq:U-x-1}
\end{align}

We note that
\[
\mathcal{A}U(Z_{s})+U'(Z_{s})(u_{i}(t,Z_{s})+u_{j}(Z_{s}))\le k_{i}u_{i}(t,Z_{s})-\pi_{i}(Z_{s})
\]
 because of condition (iii). Regarding player $j$'s singular control,
we have $U(Z_{s^{-}}+\Delta\xi_{js})-U(Z_{s^{-}})=0$ for any $\Delta\xi_{js}>0$
because of condition (i) that $U'(x)=0$ for all $x$ within the interior
of $D_{j}$. Lastly, we note that $U'(Z_{s})\sum_{l=1}^{2}d\xi_{ls}^{c}\le k_{i}d\xi_{is}^{c}$
because $U'(x)\le k_{i}$ from condition (iii) and $U'(Z_{s})d\xi_{js}^{c}=0$
from condition (i), and $U(Z_{s^{-}}+\Delta\xi_{is})-U(Z_{s^{-}})\le k_{i}\Delta\xi_{is}$
from $U'(x)\le k_{i}$. Combining all these facts, we obtain 
\begin{align*}
U(x)\ge & \mathbb{E}^{x}[e^{-r\tau_{n}\wedge t}U(Z_{\tau_{n}\wedge t})]+\mathbb{E}^{x}[\int_{0}^{\tau_{n}\wedge t}e^{-rs}(\pi_{i}(Z_{s})-k_{i}u_{i}(t,Z_{s}))ds]\\
 & -k_{i}\mathbb{E}^{x}[\int_{0}^{\tau_{n}\wedge t}e^{-rs}d\xi_{is}^{c}+\sum_{0\le s\le\tau_{n}\wedge t}e^{-rs}\Delta\xi_{is}]\;.
\end{align*}

Since $U(\cdot)$ is non-decreasing and $U'(x)\le k_{i}$, we have
$(U(Z_{t}))^{-}\le k_{i}(Z_{t})^{-}+C_{0}$ for some $C_{0}>0$. Furthermore,
in comparison to the uncontrolled process $X_{t}$, we have $Z_{t}\ge X_{t}$
for any control strategies taken by the players because the controls
$\xi_{i}$ and $\xi_{j}$ always boost $Z_{t}$. It follows that $(Z_{t})^{-}\le(X_{t})^{-}$
for all $t\ge0$ if $X_{0}=Z_{0}$. By virtue of Assumption \ref{assump:mu-sigma}
(ii), $\{e^{-r\tau}(U(Z_{\tau}))^{-}:\tau>0\;\text{stopping time}\}$
is uniformly integrable under any control strategy profile $(\xi_{i},\xi_{j})$.
Thus, we have 
\[
\liminf_{n\rightarrow\infty}\mathbb{E}^{x}[e^{-r\tau_{n}\wedge t}U(Z_{\tau_{n}\wedge t})]=\mathbb{E}^{x}[\liminf_{n\rightarrow\infty}e^{-r\tau_{n}\wedge t}U(Z_{\tau_{n}\wedge t})]=E^{x}[e^{-rt}U(Z_{t})]\:.
\]
From Fatou's lemma, we have
\begin{align*}
U(x)\ge & \liminf_{n\rightarrow\infty}\mathbb{E}^{x}[e^{-r\tau_{n}\wedge t}U(Z_{\tau_{n}\wedge t})]+\liminf_{n\rightarrow\infty}\mathbb{E}^{x}\left[\int_{0}^{\tau_{n}\wedge t}e^{-rs}(\pi_{i}(Z_{s})-k_{i}u_{i}(t,Z_{s}))ds\right]\\
 & -k_{i}\liminf_{n\rightarrow\infty}\mathbb{E}^{x}[\int_{0}^{\tau_{n}\wedge t}e^{-rs}d\xi_{is}^{c}+\sum_{0\le s\le\tau_{n}\wedge t}e^{-rs}\Delta\xi_{is}]\\
\ge & E^{x}[e^{-rt}U(Z_{t})]+\mathbb{E}^{x}[\int_{0}^{t}e^{-rs}(\pi_{i}(Z_{s})-k_{i}u_{i}(t,Z_{s}))ds]\\
 & -k_{i}\mathbb{E}^{x}[\int_{0}^{t}e^{-rs}d\xi_{is}^{c}+\sum_{0\le s\le t}e^{-rs}\Delta\xi_{is}]\;.
\end{align*}

We note that $U(Z_{t})\le M$ for some $M>0$ because it is bounded
from above. Furthermore, because $U'(\cdot)$ is bounded from above
by $k_{i}$ by condition (iii), we have $(U(Z_{t}))^{-}\le k_{i}(Z_{t})^{-}+C_{0}$
for some constant $C_{0}>0$. Recall that $(Z_{t})^{-}\le(X_{t})^{-}$
for the uncontrolled process $X_{t}$ and that $\lim_{t\rightarrow\infty}\mathbb{E}^{x}[e^{-rt}(X_{t})^{-}]=0$
by Assumption \ref{assump:mu-sigma}. Therefore, $\lim_{t\rightarrow\infty}E^{x}[e^{-rt}U(Z_{t})]=0$
is also satisfied, and so we obtain
\begin{align*}
U(x) & \ge\mathbb{E}^{x}[\int_{0}^{\infty}e^{-rs}(\pi_{i}(Z_{s})-k_{i}u_{i}(t,Z_{s}))ds]-k_{i}\int_{0}^{\infty}e^{-rs}d\xi_{is}^{c}-k_{i}\sum_{0\le s<\infty}e^{-rs}\Delta\xi_{is}]\\
 & =V_{i}(x;\xi_{i},\xi_{j})\quad.
\end{align*}
 Since $\xi_{i}$ is an arbitrary strategy of player $i$ that satisfies
$d\xi_{it}=u_{i}(t,Z_{t})+d\xi_{it}^{c}+\Delta\xi_{it}$, we have
proved that $U(x)$ dominates all payoff functions that belong to
the set $\{V_{i}(x;\xi_{i},\xi_{j}):\xi_{i}\in\Sigma_{i}\}$.

Lastly, we consider $Z$ subject to the strategy profile $(\xi_{i}^{*},\xi_{j})$.
Note that $\{e^{-r\tau}U(Z_{\tau}):\tau>0\;\text{stopping time}\}$
is uniformly integrable because $U(\cdot)$ is bounded from above
and $\{e^{-r\tau}(U(Z_{\tau}))^{-}:\tau>0\;\text{stopping time}\}$
is uniformly integrable. By condition (iv), it is straightforward
to verify that all the weak inequalities above can be exactly replaced
by equalities if $\xi_{i}$ above is replaced by $\xi_{i}^{*}$. Therefore,
$\xi_{i}^{*}$ is the best response among $\Sigma_{i}$ against $\xi_{j}$,
and $U(\cdot)$ is the best payoff function of player $i$ within
the set $\{V_{i}(x;\xi_{i},\xi_{j}):\xi_{i}\in\Sigma_{i}\}$. \eproof

\textbf{Proof of Theorem \ref{thm:symm-eq}}: As a first step, we
verify that $V(x)$ in (\ref{eq:Vx}) and $\xi$ given by the proposition
satisfy the conditions (i)\textendash (iv) of Theorem \ref{thm:verif}.

(i) Note that $C_{i}=C_{j}=I$, and $V(\cdot)\in C^{2}(I)$ and that
$V(\cdot)$ is non-decreasing. Thus, condition (i) is satisfied. Condition
(ii) is not applicable because $D_{1}=D_{2}=\emptyset$.

(iii) and (iv) For $x>\theta$, we have $V'(x)\le k$ by the definition
of $V(\cdot)$, and $\mathcal{A}V(x)+\pi(x)=0$, so condition (iii)
is satisfied. For $x\le\theta$, we have $V'(x)=k$ and
\[
\mathcal{A}V(x)+\pi(x)+u(x)V'(x)+v(U'(x)-k)=0
\]
for any arbitrary $v$ by the definition of $u(\cdot)$. Thus, condition
(iii) is also satisfied for $x\le\theta$. Because $D_{1}=D_{2}=\emptyset$,
it also follows that condition (iv) is satisfied. 

Thus, we have proved that $\xi=(\xi_{i},\xi_{j})$ is a strategy profile
that belongs to $\Sigma=\Sigma_{i}\times\Sigma_{j}$ and that $V_{i}(x;\xi)=\sup_{\zeta_{i}\in\Sigma_{i}}V_{i}(x;\zeta_{i},\xi_{j})=V(x)$
for both $i=1,2$. This does not necessarily mean that $\xi$ is a
Nash equilibrium. More precisely, we need to prove that $\xi_{i}$
is the best response among $\Xi_{i}$. To prove it, we only need to
show that $V(x)=\sup_{\zeta_{i}\in\Xi_{i}}V_{i}(x;\zeta_{i},\xi_{j})$,
i.e., that $V(\cdot)$ is player $i$'s optimal value function given
$\xi_{j}$. We do so by showing that there exists a singular control
strategy $\xi_{i}^{*}\in\Xi_{i}$ such that $V(x)=V_{i}(x;\xi_{i}^{*},\xi_{j})$
and that $V(\cdot)$ satisfies the optimality conditions for player
$i$ given in Section \ref{subsec:Single}.

Given $\xi_{j}$, the SDE of $Z$ and its $r$-excessive characteristic
operator are given by 
\begin{align*}
dZ_{t} & =[\mu(Z_{t})+u_{j}(Z_{t})]dt+\sigma(Z_{t})dW_{t}\:,\\
\mathcal{A}_{\xi_{j}} & =\frac{1}{2}\sigma(x)^{2}\partial_{x}^{2}+[\mu(x)+u_{j}(x)]\partial_{x}-r\:.
\end{align*}
Consider $\xi_{i}^{*}$ with a singular control region $D_{i}=(a,\theta]$,
which is consistent with the fact that $V'(x)=k$ for all $x\in D_{i}$.
First, note that $\mathcal{A}_{\xi_{j}}V(x)+\pi(x)=0$ and $V'(x)\le k$
for all $x\in I$. Second, it follows that $[\mathcal{A}_{\xi_{j}}V(x)+\pi(x)][V'(x)-k]=0$
for all $x\in I$. Thus, all the conditions of the optimality are
satisfied, and we conclude that $V(x)=\sup_{\zeta_{i}\in\Xi_{i}}V_{i}(x;\zeta_{i},\xi_{j})$.
\eproof

\textbf{Proof of Theorem \ref{thm:asymmetric}}: Assume that there
exists a regular control MPE with a payoff function $V_{i}(\cdot;\xi)\in C^{2}(I)$.
To prove Theorem \ref{thm:asymmetric}, we establish the following
two statements: (i) The control regions of both players must coincide.
(ii) Player $i$'s control region must be $(a,\theta_{i})$. Given
(i) and (ii), because of the assumption $\theta_{i}\neq\theta_{j}$,
we arrive at a contradiction and hence prove the theorem.

Below we employ Theorem 11.2.1 of \citet{Oksendal2003} to prove (i)
and (ii). For the theorem to be applicable, a few conditions have
to be satisfied. First, $V_{i}(\cdot;\xi)\in C^{2}(I)$ has to be
satisfied, which is assume above. Second, we need to have $\vert\mathcal{A}V_{i}(x;\xi)\vert<\infty$
$\forall x\in I$, which is satisfied because of Assumption \ref{assump:mu-sigma}.
Lastly, we need to have $\vert[u_{i}(x)+u_{j}(x)]V_{i}'(x;\xi)\vert<\infty$
$\forall x\in I$, which is satisfied because we limit $\xi$ to $\Sigma_{i}\times\Sigma_{j}$.
Thus, Theorem 11.2.1 of \citet{Oksendal2003} is applicable. 

(i) We first prove that the \emph{regular control regions} $E_{i}=\{x\in I:u_{i}(x)>0\}$
of both players must coincide, i.e., $E_{1}=E_{2}$. Suppose that
there exists a non-empty open set $F_{i}\subset E_{i}\backslash E_{j}$
such that $u_{i}(x)>0$ but $u_{j}(x)=0$ whenever $x\in F_{i}$.
Note that $V_{i}'(x;\xi)\le k_{i}$; if there exists an interval in
which $V_{i}'(x;\xi)>k_{i}$, then it behooves player $i$ to adopt
a singular control strategy in this interval, which contradicts the
assumption that the equilibrium is characterized only by regular control
strategies. By Theorem 11.2.1 of \citet{Oksendal2003}, $V_{i}(x;\xi)$
must satisfy the following HJB equation in $F_{i}$:
\[
\mathcal{A}V_{i}(x;\xi)+\pi_{i}(x)+u_{i}(x)[V_{i}'(x;\xi)-k_{i}]=0
\]
such that $u_{i}(x)>0$ only if $V_{i}'(x;\xi)=k_{i}$. Here we used
the fact that $u_{j}(x)=0$ in $F_{i}$. By the assumption that $u_{i}(x)>0$
for all $x\in F_{i}$, $V_{i}(x;\xi)=v_{0}+k_{i}x$ for some constant
$v_{0}$. Then the solution to 
\[
\mathcal{A}V_{i}(x;\xi)+\pi_{i}(x)=\mu(x)k_{i}-r[v_{0}+k_{i}x]+\pi_{i}(x)=0
\]
 cannot be a non-empty open interval according to Assumption \ref{assum:pi-p-rk}.
Thus, the necessary HJB condition cannot be satisfied in $F_{i}$.
It follows that a non-empty open set $F_{i}\subset E_{i}\backslash E_{j}$
cannot exist for either $i$. Because $u_{i}(\cdot)$ is Lipschitz
continuous, it implies $E_{i}=E_{j}$. For convenience, we let $E=E_{i}=E_{j}$
denote the common regular control region for the remainder of this
proof.

(ii) By Theorem 11.2.1 of \citet{Oksendal2003}, $V_{i}(\cdot;\xi)$
must satisfy 
\begin{equation}
\mathcal{A}V_{i}(z;\xi)+\pi_{i}(z)+u_{j}(z)\partial_{z}V_{i}(z;\xi)+u_{i}(z)[\partial_{z}V_{i}(z;\xi)-k_{i}]=0\:,\label{eq:HJB-proof}
\end{equation}
where $u_{i}(x)>0$, $V_{i}'(x;\xi)=k_{i}$ and $u_{j}(x)>0$ must
be satisfied only if $x\in E$, and $u_{i}(x)=u_{j}(x)=0$ and $V_{i}'(x;\xi)<k_{i}$
must be satisfied for $x\not\in E$. Furthermore, since $u_{j}(x)>0$
for $x\in E$, in order for (\ref{eq:HJB-proof}) to hold, $\mathcal{A}V_{i}(x;\xi)+\pi_{i}(x)<0$
needs to be satisfied for $x\in E$. In summary, the most salient
necessary conditions are $[\mathcal{A}V_{i}(x;\xi)+\pi_{i}(x)][V_{i}'(x;\xi)-k_{i}]=0$,
$\mathcal{A}V_{i}(x;\xi)+\pi_{i}(x)\le0$, and $V_{i}'(x;\xi)\le k_{i}$.
These conditions exactly coincide with the optimality conditions for
a single decision maker singular stochastic control problem given
in Section \ref{subsec:Single}. By virtue of Lemma \ref{lemma:OptSol},
there is a unique function $V_{i}^{*}(\cdot)$ given by (\ref{eq:Vx})
where $\theta$ and $\pi$ are replaced by $\theta_{1}$ and $\pi_{1}$
that satisfies these necessary conditions. Based on the form of $V_{i}^{*}(\cdot)$
given by (\ref{eq:Vx}), the regular control region is $E_{i}=(a,\theta_{i})$. 

From (ii), we conclude that the equilibrium is characterized by the
player's control region $E_{i}=(a,\theta_{i})$ and $E_{j}=(a,\theta_{j})$.
However, as established by (i), the two control regions must coincide
($E_{i}=E_{j}$), which is not possible if $\theta_{i}\neq\theta_{j}$.
Therefore, there is no regular control strategy equilibrium with associated
payoff functions $V_{i}(\cdot;\xi)\in C^{2}(I)$. \eproof 

\textbf{Proof of Proposition \ref{prop:asymm-eq}}: Note that $\xi_{1}$
continuously evolves in $C_{1}=(a,\theta')\cup(\theta,b)$ and $\xi_{2}$
continuously evolves in $C_{2}=I$ according to the strategy profile
$\xi$. To prove the proposition, it is sufficient to verify that
$U_{1}(\cdot)$ and $U_{2}(\cdot)$ with the strategy profile $\xi$
satisfy all the conditions of Theorem \ref{thm:verif}. The complete
proof exactly parallels that of Theorem \ref{thm:symm-eq}.

We first assume $\xi_{2}$ and examine $U_{1}(\cdot)$. First, it
is straightforward to verify that $U_{1}(\cdot)$ satisfies conditions
(i) and (ii) of Theorem \ref{thm:verif} because $C_{2}=I$. Hence,
we only need to verify (iii) and (iv). Because of the forms of $u_{1}(\cdot)$
and $u_{2}(\cdot)$ given in (\ref{eq:u1-asymm}) and (\ref{eq:u2-asymm}),
we have $\mathcal{A}U_{1}(x)+\pi_{1}(x)+u_{2}(x)U_{1}'(x)=0$ for
all $x<\theta'$. We also have $U_{1}'(x)\le k_{1}$ for all and $x\in I$,
and, in particular, $U_{1}'(x)=k_{1}$ for all $x<\theta_{1}$. Hence,
conditions (iii) and (iv) are satisfied. 

Next, we verify that $U_{2}(\cdot)$ given $\xi_{1}$ satisfies the
conditions (i)-(iv) of Theorem \ref{thm:verif}.

(i) Note that $D_{1}=[\theta',\theta_{1}]$ and that $U_{2}(\cdot)$
is constant in $D_{1}$. Furthermore, by the definition of $B$ in
(\ref{eq:B-theta}), we have $\lim_{x\downarrow\theta_{1}}U_{2}'(x)=0$.
Note that, by the nature of a singular control, $\xi_{1t}^{l}$ evolves
only at the right-most boundary of the singular control region $D_{1}$.
This implies that $U_{2}'(Z_{t})=0$ whenever $d\xi_{1t}^{l}>0$ so
that $U_{2}'(Z_{t})d\xi_{1t}^{l}=0$. Hence, (i) is satisfied.

(ii) The first derivative of $U_{2}(\cdot)$ is discontinuous (not
defined) at $x=\theta'$, and its second derivative is in general
not defined at $x=\theta_{1}$. However, given the first and second
derivatives of $U_{2}(\cdot)$ near $\theta'$ and $\theta_{1}$,
it is always possible to construct a function $\tilde{U}_{2}(\cdot)\in C^{2}(I)$
such that $\tilde{U}_{2}(x)=U(x)$ for all $x\in C_{1}=(a,\theta')\cup(\theta_{1},b)$
as long as $\theta'<\theta_{1}$. 

(iii) For $x>\theta_{1}$, we have $\mathcal{A}U_{2}(x)+\pi_{2}(x)=0$
by the form of $U_{2}(\cdot)$, and $U_{2}'(x)\le k_{2}$ is satisfied
so that $v_{2}[U_{2}'(x)-k_{2}]\le0$ for any $v_{2}\ge0$. Hence,
condition (iii) is satisfied for $(\theta_{1},b)$. For all $x<\theta'$,
we have $\mathcal{A}U_{2}(x)+\pi_{2}(x)+u_{1}(x)k_{2}=0$ by the form
of $u_{1}(\cdot)$, and $U_{2}'(x)=k_{2}$. Thus, condition (iii)
is satisfied for all $x\in C_{1}$. Furthermore, it is straightforward
to verify (iv) from the forms of $u_{i}(\cdot)$. \eproof

\section{R\&D Game with Spillover\label{sec:R=000026D-Game}}

An example of variable contribution games is an R\&D game with spillover,
which often occurs in high-tech industries. The technological advances
made by one firm most often spill over to another firm through various
means such as reverse engineering, leakage of information due to geographic
proximity, etc. 

We consider two firms engaging in R\&D to develop a new technology.
For simplicity, we assume that the outcome of one firm's successful
completion of R\&D completely spills over to the other.\footnote{The assumption of complete spillover is not an essential one; partial
spillover can be easily modeled, but it would complicate the analysis
without altering the main insight. } We also assume that the two firms are not in direct competition with
each other because they are in two separate markets although they
use the same technology. This model is an extension of an attrition
game: each firm would rather that its opponent conducts the R\&D.
Thus, the R\&D effort is subject to a free rider problem. Unlike the
canonical attrition game, the two firms' levels of R\&D effort are
the state variables.

Let $R$ denote the reward to each firm from the new technology, irrespective
of which firm develops it, and let $\lambda_{it}\ge0$ denote the
effort level of firm $i\in\{1,2\}$ at time $t$. The completion time
of firm $i$'s R\&D is an exponential random variable with the instantaneous
arrival rate of $\lambda_{it}$. Hence, by the property of a Poisson
process, the instantaneous arrival rate of the first completion of
R\&D is given by $\lambda_{1t}+\lambda_{2t}$. We model the cost of
maintaining the effort level of $\lambda_{it}$ as $c\lambda_{it}^{2}/2$
per unit time. Furthermore, in order to increase the effort level
by $\Delta\lambda_{it}>0$, firm $i$ has to spend $k\Delta\lambda_{it}$
for some $k>0$. Each firm $i$ can increase $\lambda_{it}$ by any
amount at any time but can never decrease it.

Let $(\lambda_{it})_{t\ge0}$ denote the non-decreasing process of
the effort level of firm $i$, and let $\hat{T}_{i}$ denote the random
completion time of firm $i$'s R\&D given the process $(\lambda_{it})_{t\ge0}$.
Given the strategy profile $\Lambda=((\lambda_{it})_{t\ge0},(\lambda_{jt})_{t\ge0})$,
firm $i$'s payoff at $t=0$ is given by 
\begin{equation}
V_{i}(\Lambda)=\mathbb{E}\left[\int_{0}^{\hat{T}_{1}\wedge\hat{T}_{2}}e^{-rt}(-\frac{c}{2}\lambda_{it}^{2})dt+e^{-r\hat{T}_{1}\wedge\hat{T}_{2}}R-\int_{0}^{\hat{T}_{1}\wedge\hat{T}_{2}}e^{-rt}kd\lambda_{it}\right]\:.\label{eq:V-Lm}
\end{equation}

The goal of this section is to demonstrate the existence of a subgame
perfect equilibrium characterized by singular control strategies.
Hence, we focus on a symmetric equilibrium in which the firms immediately
boost $\lambda_{it}$ to the equilibrium level. In particular, we
show that there exists a symmetric subgame perfect equilibrium in
which both firms immediately set their effort level at a unique value
$\lambda^{*}$ given by
\begin{equation}
\lambda^{*}=\frac{r}{8k+3c}\left[-(4k+c)r+\sqrt{(4k+c)^{2}+2(R/r-k)(8k+3c)}\right]\:,\label{eq:lm*}
\end{equation}
and maintain it until the end of the game. We can verify that this
is an equilibrium by the first-order optimality condition for each
player's best response. Here we assume that the reward is sufficiently
large so that 
\[
(4k+c)^{2}+2(R/r-k)(8k+3c)>0\:.
\]

\begin{proposition} Suppose that the initial effort levels are given
by $\lambda_{i0}<\lambda^{*}$ and $\lambda_{j0}<\lambda^{*}$. Then
there exists a subgame perfect equilibrium in which both firms boost
their effort levels up to $\lambda^{*}$ at time zero. \end{proposition}

\textbf{Proof}: Suppose that firm 2's strategy is to boost the effort
level to a level $\lambda_{2}$ at time zero and keep it at this level
until the end of the game. Our goal is to obtain the best response
of firm 1. Let $\lambda_{1,0}$ denote the initial effort level of
firm 1. Given firm 2's strategy, firm 1's best response should be
to similarly boost the effort level to some value $\lambda$ and keep
it until the end of the game because of the Markov property of the
Poisson process.

As a first step, we compute the following:
\begin{align*}
\mathbb{E}[\exp(-r\hat{T}_{1}\wedge\hat{T}_{2})]= & \int_{0}^{\infty}\lambda_{2}e^{-\lambda_{2}t_{2}}dt_{2}\int_{0}^{t_{2}}\lambda e^{-\lambda t_{1}}dt_{1}e^{-rt_{1}}\\
 & +\int_{0}^{\infty}\lambda e^{-\lambda t_{1}}dt_{1}\int_{0}^{t_{1}}\lambda_{2}e^{-\lambda_{2}t_{2}}dt_{2}e^{-rt_{2}}\\
= & \frac{\lambda+\lambda_{2}}{r+\lambda+\lambda_{2}}\:.
\end{align*}
From () we obtain
\[
V_{1}(\Lambda)=-\frac{c}{2r}\lambda^{2}+(R+\frac{c}{2r}\lambda^{2})\frac{\lambda+\lambda_{2}}{r+\lambda+\lambda_{2}}-k(\lambda-\lambda_{1,0})\:.
\]
 It follows that the first derivative is
\begin{align*}
\frac{\partial V_{1}(\Lambda)}{\partial\lambda} & =\frac{\frac{c}{2}\lambda^{2}+Rr-c\lambda(r+\lambda+\lambda_{2})}{(r+\lambda+\lambda_{2})^{2}}-k\\
 & =\frac{Rr-c\lambda(r+\lambda/2+\lambda_{2})-k(r+\lambda+\lambda_{2})^{2}}{(r+\lambda+\lambda_{2})^{2}}\:.
\end{align*}
Note that the numerator of the second line is a concave function of
$\lambda$, so the maximum value of $V_{1}(\Lambda)$ is achieved
by a unique value of $\lambda$ that satisfies the first-order condition
$\partial V_{1}(\Lambda)/\partial\lambda=0$. Assuming a symmetric
equilibrium with $\lambda=\lambda_{2}$ that solves the first order
equation $\partial V_{i}(\Lambda)/\partial\lambda=0$, we obtain $\lambda=\lambda^{*}$
given by (\ref{eq:lm*}). Therefore, immediately boosting the effort
level to $\lambda^{*}$ is the best response to firm 2's strategy
of $\lambda_{2}=\lambda^{*}$. Since the firms are symmetric, the
same is true for firm 2. We conclude that the given strategy profile
is a subgame perfect equilibrium. \eproof

Even though the game is an extension of an attrition game, the equilibrium
obtained above is characterized by an immediate lump sum (singular)
control rather than a mutual delay of action or gradualism. Intuitively,
the emergence of a singular control equilibrium is due to the players'
ability to control its states (level of effort). In contrast to canonical
attrition games that allow for only binary actions, the players of
a variable concession game can control their degree of concession
immediately by a modest amount, so they do not need to delay their
concession.


\end{document}